\documentclass[10pt]{article}
\usepackage{latexsym}
\usepackage[all]{xy}
\usepackage{amsmath}
\usepackage{amssymb}
\usepackage{amsfonts}

\newtheorem{thm}{Th\'eor\`eme}[section]
\newtheorem{prop}[thm]{Proposition}
\newtheorem{lem}[thm]{Lemme}
\newtheorem{rmk}[thm]{Remarque}

\newtheorem{df}[thm]{D\'efinition}
\newtheorem{cor}[thm]{Corollaire}

\begin{document}

\title{\textbf{Descente fid\`element plate pour les n-champs d'Artin}}
\bigskip
\bigskip

\author{\bigskip\\
Bertrand To\"en \\
\small{Institut de Math\'ematiques et de Mod\'elisation de Montpellier, UMR CNRS 5149}\\
\small{Universit\'e de Montpellier 2}\\
\small{Case Courrier 051}\\
\small{Place Eug\`ene Bataillon}\\
\small{34095 Montpellier Cedex}\\
\small{France}\\
\small{e-mail: btoen@math.univ-montp2.fr}}

\date{D\'ecembre 2009}

\maketitle

\begin{abstract}
Nous montrons deux r\'esultats de descente fid\`element plate de pr\'esentation finie dans le cadre des $n$-champs d'Artin.
Tout d'abord, un champ pour la topologie \'etale et qui est un $n$-champ d'Artin au sens de \cite{si,hagII} est aussi un champ pour la topologie fppf. De plus, un $n$-champ pour la topologie fppf et qui poss\`ede un $n$-atlas fppf est un $n$-champ d'Artin (i.e. poss\`ede aussi un $n$-atlas lisse). Nous d\'eduisons de ces deux r\'esultats un th\'eor\`eme de comparaison entre cohomologies \'etale et fppf (\`a coefficients dans des sch\'emas en groupes non n\'ecessairement
lisses ou encore non-ab\'eliennes). Ce travail est \'ecrit dans le contexte des champs d\'eriv\'es de \cite{hagII,seat}, et ces r\'esultats vallent donc aussi pour des $n$-champs d'Artin d\'eriv\'es.
\end{abstract}

\section*{Introduction}

L'objectif de ce travail est de d\'emontrer l'\'equivalence de deux notions, \`a priori diff\'erentes, de 
$n$-champs alg\'ebriques\footnote{Par la suite, l'expression \emph{champs} fera toujours r\'ef\'erence \`a la notion de champs sup\'erieurs. Les champs en groupo\"\i des seront alors appel\'es des \emph{$1$-champs}.}. Comme cela est expliqu\'e dans \cite[\S 1.3]{hagII}, il existe une notion 
g\'en\'erale de \emph{champs ($n$-)g\'eom\'etriques}, qui d\'epend du choix d'un couple
$(\tau,\mathbf{P})$, form\'e d'une topologie de Grothendieck $\tau$ sur le site des sch\'emas affines, et d'une classe de morphismes $\mathbf{P}$ entre sch\'emas affines, et satisfaisant \`a certaines conditions de compatibilit\'e. Par d\'efinition, les champs g\'eom\'etriques pour le coupe $(\tau,\mathbf{P})$
sont obtenus en recollant des sch\'emas affines le long de relations d'\'equivalences it\'er\'ees de 
type $\mathbf{P}$. Si $\tau=et$ est la topologie \'etale, et 
$\mathbf{P}=et$ la classe des morphismes \'etales, la notion de champs g\'eom\'etriques correspondante est 
celle de champs alg\'ebriques de Deligne-Mumford, ou plus pr\'ecis\'ement son extension au cadre des 
champs sup\'erieurs. Si $\tau=et$ est la topologie \'etale et $\mathbf{P}=li$ est la classe des morphismes
lisses, alors la notion de champs g\'eom\'etriques correspond \`a celle de 
champs alg\'ebriques d'Artin (au sens de \cite{si} et de \cite[\S 2.1]{hagII}). Dans ce travail nous nous int\'eresserons au cas du couple $(fppf,pl)$, form\'e de la topologie fid\`element plate et de pr\'esentation finie (not\'ee $fppf$), et de la classe
$pl$ des morphismes plats de pr\'esentation finie. Notre r\'esultat principal affirme que les champs
g\'eom\'etriques pour le couple $(fppf,pl)$ sont exactement les champs g\'eom\'etriques
pour le couple $(et,li)$, c'est \`a dire les champs alg\'ebriques d'Artin.

\begin{thm}\label{ti}
Le foncteur \emph{champ associ\'e pour la topologie fppf} induit une \'equivalence de la cat\'egorie
des champs g\'eom\'etriques pour le couple $(et,li)$ avec celle des champs g\'eom\'etriques
pour le couple $(fppf,pl)$.
\end{thm}

De mani\`ere \'equivalente, le th\'eor\`eme pr\'ec\'edent peut aussi s'\'enoncer en deux 
assertions, l'une concernant la pleine fid\`elit\'e et l'autre l'essentielle surjectivit\'e:

\begin{enumerate}
\item Un champ g\'eom\'etrique pour le couple $(et,li)$ est un champ pour la topologie $fppf$.
\item Un champ pour la topologie $fppf$, qui poss\`ede un atlas plat et localement de pr\'esentation finie, poss\`ede un atlas lisse. 
\end{enumerate}

La premi\`ere assertion affirme que la cohomologie non-ab\'elienne, pour la topologie \'etale et \`a valeurs dans un 
champs alg\'ebrique d'Artin, peut aussi se calculer en utilisant  la topologie $fppf$. Il s'agit donc d'un \'enonc\'e de comparaison entre cohomologie \'etale et $fppf$, qui g\'en\'eralise le fait bien connu $H^{i}_{et}(X,G)\simeq H^{i}_{fppf}(X,G)$, pour $X$ un sch\'ema et $G$ un sch\'ema en groupes ab\'eliens lisse (voir
\cite[App. 11]{gb}). La seconde assertion est une g\'en\'eralisation au cadre des champs sup\'erieurs d'un th\'eor\`eme d'Artin qui affirme que le $1$-champ quotient d'un groupo\"\i de plat et de pr\'esentation finie
dans les sch\'emas peut aussi s'\'ecrire comme le $1$-champ quotient d'un groupo\"\i de lisse (voir par exemple
\cite{lm}). Comme nous le verrons, il se trouve que le point $(1)$ est une cons\'equence du point $(2)$
(voir notre lemme \ref{l1}). Le contenu du th\'eor\`eme \ref{ti} est donc essentiellement le fait que l'existence d'un atlas plat et localement de pr\'esentation finie implique l'existence d'un atlas lisse. 

Le th\'eor\`eme \ref{t1} poss\`ede plusieurs cons\'equences int\'eressantes. L'une est l'existence d'une 
th\'eorie des \emph{$n$-gerbes alg\'ebriques}, g\'en\'eralisation d'ordre sup\'erieur des 
gerbes alg\'ebriques de \cite{lm}, et de l'existence de gerbes r\'esiduelles (ce qui implique, entre autre, la repr\'esentabilit\'e des faisceaux d'homotopie). Nous ne pr\'esenterons pas 
ces cons\'equences dans ce travail, et nous renvoyons le lecteur \`a \cite[\S 2.2, \S 2.3]{kch}, o\`u il trouvera, de plus, des applications aux invariants cohomologiques des $n$-champs d'Artin. En contre partie, nous avons inclu des applications \`a la comparaison entre cohomologies (non-ab\'eliennes) \'etales et $fppf$. Pour 
$G$ un champ en groupes lisses, nous montrons que les espaces $H^{1}_{et}(X,G)$ et $H^{1}_{fppf}(X,G)$
coincident, g\'en\'eralisant ainsi l'\'enonc\'e analogue pour des sch\'emas en groupes. Travailler avec des champs sup\'erieurs nous permet d'introduire la notion de champs en $k$-groupes $G$, qui poss\`ede des champs
classifiants $K(G,i)$ pour tout $i\leq k$, ce qui permet ainsi de d\'efinir
les espaces $H^{i}_{\tau}(X,G)$ pour tout $i\leq k$ (pour $\tau$ une des deux topologies et ou fppf). Lorsqu'un tel $G$ est plat de localement de pr\'esentation finie le th\'eor\`eme \ref{ti} implique que $K(G,i)$ est un champ g\'eom\'etrique, lisse lorsque 
$i>0$. Nous tirons de cela, pour tout sch\'ema en groupes ab\'eliens plats et localement de pr\'esentation finie $G$, et tout sch\'ema $X$, l'existence d'une suite exacte longue
$$\xymatrix{
H^{i-2}_{et}(X,\underline{H}^{1}_{fppf}(-,G)) \ar[r] & H^{i}_{et}(X,G) \ar[r] & H^{i}_{fppf}(X,G) \ar[r] & H^{i-1}_{et}(X,\underline{H}^{1}_{fppf}(-,G)) \ar[r] & 
H^{i+1}_{et}(X,G), }$$
o\`u $\underline{H}^{1}_{fppf}(-,G)$ est le faisceau, pour la topologie \'etale, associ\'e 
au pr\'efaisceau $U \mapsto H^{1}_{fppf}(U,G)$. Ceci s'exprime aussi en disant que 
$\mathbb{R}^{i}f_{*}(G)=0$ pour tout $i>1$, o\`u $f$ est le morphisme g\'eom\'etrique
de passage de la topologie $fppf$ \`a la topologie \'etale. Dans le m\^eme genre d'id\'ee, si 
$A$ est un anneau hens\'elien excellent de corps r\'esiduel $k$, et si 
$G$ est un sch\'ema en groupes plat et localement de pr\'esentation finie sur $A$, de fibre
sp\'eciale $G_{0}$, alors nous montrons que le morphisme de restriction
$$H^{i}_{fppf}(Spec\, A,G) \longrightarrow H^{i}_{fppf}(Spec\, k,G_{0})$$
est un isomorphisme pour $i>1$, et surjectif pour $i=1$. 
Nous n'avons pas trouv\'e d'\'enonc\'es dans ce go\^ut dans la litt\'erature, et il est possible que ces deux 
r\'esultats soit nouveaux. \\

La strat\'egie de preuve du th\'eor\`eme \ref{ti} est tout \`a fait analogue \`a celle utilis\'ee
dans \cite{lm} pour traiter le cas des $1$-champs. Le point cl\'e est la construction, \`a partir d'un 
recouvrement plat et localement de pr\'esentation finie entre sch\'emas $p : X' \rightarrow X$, d'un 
recouvrement lisse $Y \rightarrow X$. Le sch\'ema $Y$ est le sch\'ema des \emph{quasi-sections quasi-lisses} de la projection $p$, qui classifie les extensions finies plates $Z \rightarrow X$ munies d'un morphisme
l.c.i. $Z \rightarrow X'$ au-dessus de $X$. Le principal apport de ce travail est de montrer que cette construction 
reste raisonnable lorsque $p$ est 
maintenant un atlas $fppf$ de $n$-champs, ce qui complique consid\'erablement les d\'etails techniques. Cependant, les \'etapes de la preuve que nous donnons suivent essentiellement celles dans le cas des sch\'emas, avec m\^eme une simplification pour d\'emontrer la lissit\'e de $Y$, due \`a l'utilisation
du langage de la g\'eom\'etrie alg\'ebrique d\'eriv\'ee (ce qui permet de ramener la preuve de la lissit\'e 
de $Y \rightarrow X$
au \emph{simple calcul} d'un espace cotangent). L'int\'egralit\'e de ce travail est d'ailleur \'ecrit dans
le langage des champs d\'eriv\'es de \cite{hagII}, et notre th\'eor\`eme \ref{ti} est donc aussi valable
dans ce contexte. \\

\bigskip
\textbf{Remerciements:} Le r\'esultat principal de ce travail avait \'et\'e annonc\'e, avec quelques rudiments de preuve, dans \cite{seat}. C'est la remarque, raviv\'ee r\'ecemment par une correspondance avec 
Chetan Balwe, que les constructions et r\'esultats 
de \cite{kch}
en d\'ependent qui m'a insit\'e \`a en r\'ediger une preuve compl\`ete. Je remercie donc Chetan pour son message (et je ne peux qu'inviter le lecteur \`a lire sa th\`ese \cite{ba}). 

Par ailleurs, cet article a \'et\'e r\'edig\'e \`a cheval entre Toulouse et Montpellier, en grande partie dans le fameux train CORAIL TEOZ de 17h45, et je souhaite remercier la SNCF pour avoir su garder des retards tout \`a fait raisonnables durant cette p\'eriode. Ces va-et-viens font suite \`a mon changement d'affectation, de l'Institut de Math\'ematique de Toulouse vers l'Institut de Math\'ematiques et de Mod\'elisation de Montpellier: je remercie le second pour un acceuil tr\`es chaleureux. \\

\bigskip

\textbf{Notation et terminologie:} \begin{itemize}

\item $k$ : un anneau commutatif de base fix\'e, ou plus g\'en\'eralement un anneau commutatif simplicial (pour les lecteurs \emph{braves}).

\item $sk-CAlg$ : la cat\'egorie des $k$-alg\`ebres simpliciales. 

\item $dAff_{k}:=(sk-CAlg)^{op}$ : la cat\'egorie des $k$-sch\'emas affines d\'eriv\'es.

\item $\tau$ : une topologie de mod\`eles sur $dAff_{k}$, ou bien la topologie \'etale, ou bien 
la topologie $fppf$.

\item $\tau$-{champ d\'eriv\'e} : un champ sur le site de mod\`eles $(dAff_{k},\tau)$.

\item {$\tau$-champ non d\'eriv\'e (ou tronqu\'e)} : un champ sur le site 
$Aff_{k}$ des $k$-sch\'emas affines (non d\'eriv\'es) munie de la topologie $\tau$ induite.

\item $dSt_{\tau}(k)$ : la cat\'egorie des $\tau$-champs d\'eriv\'es. 

\item $St_{\tau}(k)$ : la cat\'egorie des $\tau$-champs non d\'eriv\'es, vue comme sous-cat\'egorie pleine 
de $dSt_{\tau}(k)$.

\end{itemize}

\section{Deux notions de champs d\'eriv\'es $n$-g\'eom\'etriques}

Dans cette section nous rappelons la notion de $n$-champs g\'eom\'etriques introduite dans
\cite[\S 1.3]{hagII}. Nous travaillerons au-dessus de la cat\'egorie de mod\`eles des $k$-alg\`ebres simpliciales, dont nous ferons varier la topologie $\tau$ ainsi que la classe de morphismes $\mathbf{P}$ utilis\'ee pour d\'efinir les atlas. Le lecteur trouvera le formalisme g\'en\'eral des champs g\'eom\'etriques au-dessus d'un HAG contexte dans \cite{hagII}, dont les notions pr\'esent\'ees ici sont des 
cas particuliers. \\

\subsection{Changement de topologie et de contextes pour les champs d\'eriv\'es}

Nous notons $sk-CAlg$ la cat\'egorie des $k$-alg\`ebres simpliciales commutatives, que nous munissons de sa structure de mod\`eles standard (voir par exemple \cite[\S 2.2.1]{hagII}). Nous notons aussi $dAff_{k}:=sk-CAlg^{op}$ la cat\'egorie oppos\'ee, muni de la structure de mod\'eles induite. 
Soit $\tau$ une (pr\`e-)topologie de mod\`eles sur $dAff_{k}$ (voir \cite[\S 4.3]{hagI}). Cette  topologie donne lieu \`a une 
topologies de Grothendieck sur $Ho(dAff_{k})$ (encore not\'ees $\tau$) ainsi qu'\`a
la cat\'egorie de mod\`eles des champs sur $dAff_{k}$ qui lui correspond (voir \cite[\S 4.6]{hagI})
$dAff_{k}^{\sim,\tau}$. La cat\'egorie homotopique des champs d\'eriv\'es pour $\tau$ (nous dirons aussi 
\emph{$\tau$-champ d\'eriv\'es}, ou encore \emph{champ d\'eriv\'e} si la topologie $\tau$ est 
soit implicite, ou bien si l'\'enonc\'e ne d\'epend pas de la topologie choisie) est 
$$dSt_{\tau}(k):=Ho(dAff_{k}^{\sim,\tau}).$$
Rappelons que $dSt_{\tau}(k)$ s'identifie naturellement \`a la sous-cat\'egorie pleine de la cat\'egorie homotopique 
$Ho(SPr(dAff_{k}))$, des pr\'efaisceaux simpliciaux sur $dAff_{k}$, form\'ee des objets qui d'une part pr\'eservent les \'equivalences, et d'autre part poss\`edent la propri\'et\'e de descente pour les $\tau$-hyper-recouvrements (voir \cite[\S 4.4]{hagI}). Par la suite nous verrons toujours $dSt_{\tau}(k)$ comme plong\'ee dans 
$Ho(SPr(dAff_{k}))$.

La cat\'egorie des $\tau$-champs d\'eriv\'es contient une sous-cat\'egorie pleine
$St_{\tau}(k) \subset dSt_{\tau}(k)$, form\'ee des $\tau$-champs \emph{non-d\'eriv\'es} (nous dirons aussi 
\emph{tronqu\'es}). La cat\'egorie $St_{\tau}(k)$ est \'equivalente \`a la cat\'egorie homotopique des
pr\'efaisceaux simpliciaux sur le site des sch\'emas affines (non-d\'eriv\'es) $Aff_{k}$, muni de la
topologie $\tau$ restreinte aux $k$-alg\`ebres non-simpliciales. Le foncteur 
d'inclusion $St_{\tau}(k) \hookrightarrow  dSt_{\tau}(k)$ est alors obtenu par 
extension de Kan \`a gauche des pr\'efaisceaux simpliciaux le long de l'inclusion 
naturelle $Aff_{k} \hookrightarrow dAff_{k}$ induite par l'inclusion des $k$-alg\`ebres dans les 
$k$-alg\`ebres simpliciales (qui consiste \`a voir une $k$-alg\`ebre comme une 
$k$-alg\`ebre simpliciale constante). Ce foncteur d'inclusion 
$i : St_{\tau}(k) \longrightarrow dSt_{\tau}(k)$ poss\`ede un adjoint \`a droite 
$$t_{0} : dSt_{\tau}(k) \longrightarrow St_{\tau}(k),$$
qui \`a un pr\'efaisceaux simplicial sur $dAff_{k}$ associe sa restriction \`a $Aff_{k}$.  Par la suite, 
nous identifierons syst\'ematiquement $St_{\tau}(k)$ \`a une sous-cat\'egorie pleine de $dSt_{\tau}(k)$, form\'ee des $F$ tels que le morphisme naturel $t_{0}(F) \longrightarrow F$ soit un isomorphisme. Nous renvoyons \`a \cite[\S 2.1]{hagII} pour plus de d\'etails sur les champs non-d\'eriv\'es. \\

Soient maintenant une classe de morphismes $\mathbf{P}$ dans $Ho(dAff_{k})$. Nous supposerons que la topologie $\tau$  satisfait les conditions \cite[1.3.2.2]{hagII}. Nous supposerons aussi que le couple $(\tau,\mathbf{P})$ satisfait les conditions \cite[1.3.2.11]{hagII}. Nous ne rappelons pas ces conditions ici, qui affirment, en gros, que la topologie $\tau$ est sous-canonique, compatible avec les sommes disjointes finies, et que les morphismes de $\mathbf{P}$ ont de bonnes propri\'et\'es de localit\'e par rapport \`a $\tau$. Par la suite nous nous int\'eresserons \`a deux exemples: $\tau$ sera la topologie \'etale et 
$\mathbf{P}$ les morphismes lisses, ou encore $\tau$ sera la topologie fppf et $\mathbf{P}$ sera les morphismes plats et de pr\'esentation presque finie (dont les d\'efinitions seront rappel\'ees dans les deux paragraphes suivants).

La donn\'ee du couple $(\tau,\mathbf{P})$ permet, comme cela est expliqu\'e dans \cite[Def. 1.3.3.1]{hagII}, de d\'efinir une notion de
champs d\'eriv\'es $(n,\mathbf{P})$-g\'eom\'etriques. Les champs d\'eriv\'es $(n,\mathbf{P})$-g\'eom\'etriques forment une sous-cat\'egorie pleine
de $dSt_{\tau}(k)$, not\'ee $dSt_{\tau}^{n,\mathbf{P}}(k)$. On a des inclusions naturelles
$$dSt_{\tau}^{n,\mathbf{P}}(k) \subset dSt_{\tau}^{n+1,\mathbf{P}}(k) \subset dSt_{\tau}(k),$$
et la r\'eunion de ces sous-cat\'egories sera not\'ee 
$$dSt_{\tau}^{\mathbf{P}}(k):=\cup_{n} dSt_{\tau}^{n,\mathbf{P}}(k).$$

Rappelons la d\'efinition, par induction sur $n$, des cat\'egories $dSt_{\tau}^{n,\mathbf{P}}(k)$. Pour $n=0$, la sous-cat\'egorie $dSt_{\tau}^{0,\mathbf{P}}(k)$ consiste en tous les champs d\'eriv\'es affines, c'est \`a dire de la forme 
$\mathbb{R}\underline{Spec}\, A$ pour un $A\in sk-CAlg$\footnote{Les conventions de ce travail diff\`erent de celles de \cite{hagII}, les objets affines dans \cite{hagII} \'etant $(-1)$-g\'eom\'etriques.}. Les conditions sur $\tau$ impliquent que
le foncteur 
$$\mathbb{R}\underline{Spec} : Ho(sk-CAlg)^{op} \longrightarrow  dSt_{\tau}(k)$$
est pleinement fid\`ele. Ainsi, $dSt_{\tau}^{0,\mathbf{P}}$ est naturellement \'equivalente \`a la cat\'egorie
$Ho(dAff_{k})$ (et se trouve donc \^etre ind\'ependante du choix de $\mathbf{P}$ et de $\tau$). Un morphisme entre champs
d\'eriv\'es
$f : F \longrightarrow G$ est dit $0$-g\'eom\'etrique (ou affine), si pour tout $X$ affine et tout
$X \longrightarrow G$, le champ d\'eriv\'e $F\times_{G}^{h}X$ est $0$-g\'eom\'etrique. Un tel morphisme est de plus dans $\mathbf{P}$ si la projection induite $X \longrightarrow F\times_{G}^{h}X$ correspond \`a un morphisme de $\mathbf{P}$ dans $Ho(dAff_{k})$.

Supposons maintenant $n>0$ et que la
cat\'egorie $dSt_{\tau}^{n-1,\mathbf{P}}(k)$ soit d\'efinie, ainsi que la notion de morphismes $(n-1,\mathbf{P})$-repr\'esentables et de morphismes $(n-1,\mathbf{P})$-repr\'esentables dans $\mathbf{P}$. On d\'efinit alors la sous-cat\'egorie
$dSt_{\tau}^{n,\mathbf{P}}(k)$, ainsi que les notions de morphismes $(n,\mathbf{P})$-repr\'esentables et de morphismes $(n,\mathbf{P})$-repr\'esentables et dans $\mathbf{P}$, de la fa\c{c}on suivante. 

\begin{enumerate}

\item Un champ d\'eriv\'e $F \in dSt_{\tau}(k)$ est $(n,\mathbf{P})$-g\'eom\'etrique s'il existe une famille 
$\{X_{i}\}_{i}$ d'affines, et un \'epimorphisme de champs 
$$\coprod_{i}X_{i} \longrightarrow F$$
tel que chaque morphisme $X_{i} \longrightarrow F$ soit un morphisme $(n-1,\mathbf{P})$-g\'eom\'etrique et dans $\textbf{P}$.
Une telle donn\'ee pour $F$ sera appel\'ee par la suite un \emph{$(n,\mathbf{P})$-atlas pour $F$}.

\item Un morphisme de champs d\'eriv\'es $f : F \longrightarrow G$ est $(n,\mathbf{P})$-repr\'esentable 
(nous dirons aussi $(n,\mathbf{P})$-g\'eom\'etrique) si pour tout
affine $X$, et tout morphisme $X \longrightarrow G$, le champ d\'eriv\'e
$F\times_{G}^{h}X$ est $(n,\mathbf{P})$-g\'eom\'etrique.

\item Un morphisme de champs d\'eriv\'es $f : F \longrightarrow G$ est $(n,\mathbf{P})$-repr\'esentable et 
dans \textbf{P} s'il est d'une part
$(n,\mathbf{P})$-repr\'esentable, et d'autre part si pour tout $X \longrightarrow G$ avec $X$ affine,
il existe un $(n,\mathbf{P})$-atlas
$$\coprod_{i}X_{i} \longrightarrow F\times_{G}^{h}X$$
tel que tous les morphismes induits
$X_{i} \longrightarrow X$ entre champs affines soient dans $\mathbf{P}$.

\end{enumerate}

Supposons maintenant que l'on se donne deux topologies de mod\`eles $\tau$ et $\tau'$, et deux classes de morphismes
$\mathbf{P}$ et $\mathbf{P}'$, de sorte \`a ce que les couples $(\tau,\mathbf{P})$ et $(\tau',\mathbf{P}')$ satisfassent 
tous deux aux conditions \cite[1.3.2.2,1.3.2.11]{hagII}. Nous supposerons que le couple $(\tau,\mathbf{P})$ est \emph{plus fort} que 
$(\tau',\mathbf{P}')$ au sens suivant.

\begin{enumerate}
\item Tout $\tau$-champ d\'eriv\'e est un $\tau'$-champ d\'eriv\'e.
\item On a $\mathbf{P'} \subset \mathbf{P}$.
\end{enumerate}

La condition $(1)$ ci-dessus dit qu'un pr\'efaisceau simplicial $F : dAff_{k}^{op} \longrightarrow SEns$, qui pr\'es\`erve
les \'equivalences et qui v\'erifie la condition de descente pour les $\tau$-hyper-recouvrements, v\'erifie aussi la condition
de descente pour les $\tau'$-hyper-recouvrements. En identifiant les cat\'egories $dSt_{\tau}(k)$ et 
$dSt_{\tau'}(k)$ \`a des sous-cat\'egories pleines de $Ho(SPr(dAff_{k}))$, cette condition est \'equivalente au fait que 
$dSt_{\tau}(k) \subset dSt_{\tau'}(k)$.

Consid\'erons maitenant le foncteur d'inclusion
$$i : dSt_{\tau}(k) \longrightarrow dSt_{\tau'}(k).$$
Le foncteur de champs associ\'es pour la topologie $\tau$, restreint \`a $dSt_{\tau'}(k)$, fournit un adjoint \`a gauche de ce foncteur d'inclusion
$$a : dSt_{\tau'}(k) \longrightarrow dSt_{\tau}(k).$$
Rappelons de plus que ce foncteur, ou plut\^ot son relev\'e naturel au niveau des cat\'egories de mod\`eles, commute aux limites homotopiques finies (voir \cite[Prop. 3.4.10]{hagI}). Cette propri\'et\'e d'exactitude et la condition 
$\mathbf{P'} \subset \mathbf{P}$ impliquent alors facilement (c'est \`a dire en utilisant les propri\'et\'es
\'el\'ementaires des champs g\'eom\'etrique donn\'ees dans \cite[\S 1.3.3]{hagII}), par induction sur $n$, que le 
foncteur $a$ transforme les $\tau'$-champs d\'eriv\'es $(n,\mathbf{P}')$-g\'eom\'etriques en des $\tau$-champs d\'eriv\'es
$(n,\mathbf{P})$-g\'eom\'etriques. Il induit ainsi, pour tout $n\geq 0$, un foncteur 
$$a : dSt_{\tau'}^{n,\mathbf{P}'}(k) \longrightarrow dSt_{\tau}^{n,\mathbf{P}}(k).$$

\subsection{Champs d\'eriv\'es $(n,li)$-g\'eom\'etriques et $(n,pl)$-g\'eom\'etriques}

Nous sp\'ecifions dans ce paragraphe deux couples $(\tau,\mathbf{P})$ et $(\tau',\mathbf{P}')$ comme dans le paragraphe pr\'ec\'edent. Le premier couple $(\tau,\mathbf{P})$ sera form\'e de la topologie fid\`element plate de pr\'esentation (presque) finie et des morphismes plats de pr\'esentation (presque) finie, 
et le second, $(\tau',\mathbf{P}')$, de la topologie \'etale et des morphismes 
lisses. Nous commencerons donc par rappeler les d\'efinitions des morphismes \'etales, lisses, plats et plats de pr\'esentation presque finie entre $k$-alg\`ebres simpliciales commutatives. 
Soit donc $f : A \longrightarrow B$ un morphisme dans $sk-CAlg$. On rappelle les notions suivantes, qui sont essentiellement tir\'ees de \cite{hagII} (sauf la notion $(1)$), r\'ef\'erence dans laquelle le lecteur trouvera aussi des d\'efinitions \'equivalentes en termes de complexes cotangents ou d'exactitude de foncteurs de changement de bases.

\begin{enumerate}

\item Le morphisme $f$ est \emph{de pr\'esentation presque finie} si le morphisme induit 
$\pi_{0}(A) \longrightarrow \pi_{0}(B)$ est un morphisme de pr\'esentation finie d'anneaux commutatifs.

\item Le morphisme $f$ est \emph{plat} si le morphisme $\pi_{0}(A) \longrightarrow \pi_{0}(B)$ est plat, et si de plus pour tout
$i>0$ le morphisme naturel
$$\pi_{i}(A)\otimes_{\pi_{0}(A)}\pi_{0}(B) \longrightarrow \pi_{i}(B)$$
est un isomorphisme.

\item Le morphisme $f$ est \emph{lisse} s'il est plat et si de plus $\pi_{0}(A) \longrightarrow \pi_{0}(B)$ est un morphisme lisse (en particulier de pr\'esentation finie) d'anneaux commutatifs. 

\item Le morphisme $f$ est \emph{\'etale} s'il est plat et si de plus $\pi_{0}(A) \longrightarrow \pi_{0}(B)$ est un morphisme \'etale (en particulier de pr\'esentation finie) d'anneaux commutatifs. 

\end{enumerate}

On montre que les morphismes plats, plats de pr\'esentation presque finie, lisses et \'etales sont stables par
changement de bases (homotopiques) dans $Ho(dAff_{k})$ et par composition. Cela se d\'eduit ais\'ement du fait que 
lorsque $A \longrightarrow B$ est plat, alors pour tout $A$-module simplicial $M$ le morphisme naturel
$$\pi_{*}(M)\otimes_{\pi_{0}(A)}\pi_{0}(B) \longrightarrow \pi_{*}(M\otimes_{A}^{\mathbb{L}}B)$$
est un isomorphisme.

\begin{rmk}\label{r1}
\emph{
Les morphismes lisses et \'etales d\'efinis ci-dessus sont \emph{(homotopiquement) de pr\'esentation finie} dans la cat\'egorie de mod\`eles  $sk-CAlg$ (voir \cite[Def. 1.2.3.1]{hagII}). Cependant les morphismes de pr\'esentation presque finie ne sont pas de pr\'esentation finie en ce sens, ni m\^eme les morphismes plats et de pr\'esentation presque finie. Par exemple, si $k$ est un corps, toute $k$-alg\`ebre de type finie est presque de pr\'esentation finie et plate au sens ci-dessus, mais elle n'est  homotopiquement de pr\'esentation finie dans $sk-CAlg$ que lorsqu'elle est localement d'intersection compl\`ete.}
\end{rmk}

Pour finir, une famille de morphismes $\{A \longrightarrow A_{i}\}$ dans $sk-CAlg$ est \emph{surjective}, si le morphisme de sch\'emas 
$$\{Spec\, \pi_{0}(A_{i}) \longrightarrow Spec\, \pi_{0}(A)\}$$
est surjectif (i.e. tout id\'eal premier de $\pi_{0}(A)$ se rel\`eve en un id\'eal premier d'un des $\pi_{0}(A_{i})$).
Nous dirons alors qu'une telle famille $\{A \longrightarrow A_{i}\}$ est un \emph{recouvrement plat et de pr\'esentation presque finie}, si tous les morphismes $A \longrightarrow A_{i}$ sont plats et de pr\'esentation presque finie, et si de plus la famille 
$\{A \longrightarrow A_{i}\}$ est surjective. De m\^eme, une telle famille est un \emph{recouvrement \'etale} si tous les morphismes $A \longrightarrow A_{i}$ sont \'etales, et si de plus la famille 
$\{A \longrightarrow A_{i}\}$ est surjective. On voit sans peine que les recouvrements plats et de pr\'esentation presque finie, et les recouvrement \'etale, forment deux topologie de mod\`eles sur $dAff_{k}$. Ces deux topologies seront respectivement appel\'ees les topologies \emph{plate et de pr\'esentation presque finie}, et 
\emph{\'etale}. 
La topologie plate et de pr\'esentation presque finie sera symb\^oliquement not\'ee \emph{fppf}. La topologie \'etale sera not\'ee \emph{et}. \\

Nous consid\'erons maintenant $pl$, la classe des morphismes plats et de pr\'esentation presque finie dans $dAff_{k}$, et 
$li$ celle des morphismes lisses. Les couples $(et,li)$ et $(fppf,pl)$ satisfont aux conditions de 
\cite[1.3.2.2,1.3.2.11]{hagII}. De plus, le couple $(fppf,pl)$ est clairement plus fort que le couple $(et,li)$. On dispose donc d'un foncteur entre les cat\'egories de champs d\'eriv\'es $n$-g\'eom\'etriques correspondantes
$$\phi_{n} : dSt_{et}^{n,li} \longrightarrow dSt_{fppf}^{n,pl}$$
induit par le foncteur de champ associ\'e pour la topologie $fppf$. \\

Pour terminer cette section signalons l'aspect local pour la topologie fppf de la lissit\'e. Ce r\'esultat nous sera utile par
la suite.

\begin{prop}\label{p1}
Soit $X \longrightarrow Y$ un morphisme dans $dAff_{k}$. Alors, $f$ est lisse si et seulement s'il existe un
recouvrement $fppf$ $\{X_{i}\longrightarrow X\}$ tel que les morphismes compos\'es 
$X_{i} \longrightarrow Y$ soient lisses.
\end{prop}

\textit{Preuve:} Seule la suffisance demande une preuve.
Soit $X=Spec\, B$, $Y=Spec A$, et $X_{i}=Spec\, B_{i}$. On commence par voir que 
$A \longrightarrow B$ est un morphisme plat. En effet, la famille de morphismes d'anneaux non simpliciaux $\{\pi_{0}(B) \longrightarrow \pi_{0}(B_{i})\}$
est fid\`element plate par hypoth\`ese, et les morphismes compos\'es $\pi_{0}(A) \longrightarrow \pi_{0}(B_{i})$
sont plats. Cela implique que $\pi_{0}(A) \longrightarrow \pi_{0}(B)$ est un morphisme plat. Il faut de plus montrer que 
le morphisme
$$f : \pi_{*}(A)\otimes_{\pi_{0}(A)}\pi_{0}(B) \longrightarrow \pi_{*}(B)$$
est bijectif. Mais, par changement de base le long de $\pi_{0}(B) \longrightarrow \pi_{0}(B_{i})$ ce morphisme devient
$$\pi_{*}(A)\otimes_{\pi_{0}(A)}\pi_{0}(B)\otimes_{\pi_{0}(B)}\pi_{0}(B_{i}) \simeq
\pi_{*}(A)\otimes_{\pi_{0}(A)}\pi_{0}(B_{i}) \longrightarrow \pi_{0}(B_{i}).$$
Par hypoth\`ese ces morphismes sont bijectifs, et comme $\{\pi_{0}(B) \longrightarrow \pi_{0}(B_{i})\}$
est un recouvrement fid\`element plat on en d\'eduit que $f$ est bijectif.

On vient de voir que $A\longrightarrow B$ est un morphisme plat. Il nous reste \`a montrer que
$\pi_{0}(A) \longrightarrow \pi_{0}(B)$ est aussi lisse. Pour cela on peut sans perte de g\'en\'eralit\'e supposer que 
$A=\pi_{0}(A)$ (et donc $B\simeq \pi_{0}(B)$, $B_{i}\simeq \pi_{0}(B_{i})$), et l'\'enonc\'e se ram\`ene alors \`a la localit\'e pour la topologie 
$fppf$ des morphismes lisses entre anneaux non simpliciaux. Ce dernier fait est bien 
connu (voir par exemple \cite[Prop. 17.7.7]{egaIV-4}), et peut se d\'emontrer de la fa\c{c}on suivante. Soit donc un morphisme plat $A \longrightarrow B$, et un recouvrement $fppf$ $\{B \longrightarrow B_{i}\}$, tel que chaque $A \longrightarrow B_{i}$ soit 
un morphisme lisse, tout cela pour des anneaux commutatifs non simpliciaux. 
En particulier $A \longrightarrow B_{i}$ est de pr\'esentation finie, et cela implique par 
descente fid\`element plate que $A \longrightarrow B$ est un morphisme de pr\'esentation finie. On pourra donc, par l'argument standard, se ramener au cas o\`u tous les anneaux sont de type fini sur $\mathbf{Z}$, et en particulier noeth\'eriens (voir 
par exemple \cite[Cor. 11.2.6.1]{egaIV-3}).  Il reste \`a montrer que
$A \longrightarrow B$ est aussi formellement lisse, ou de mani\`ere \'equivalente que pour tout corps 
alg\'ebriquement clos $K$ et tout morphisme $A \longrightarrow K$, la $K$-alg\`ebre $B_{K}:=B\otimes_{A}K$ est 
lisse. On se ram\`ene ainsi au cas o\`u $A=K$ est un corps alg\'ebriquement clos. On dipose ainsi 
de $B$ une $K$-alg\`ebre de type finie, d'un recouvrement $fppf$ $\{B\longrightarrow B_{i}\}$, tel que
chaque $B_{i}$ soit lisse sur $K$. Pour montrer que $B$ est lisse, il suffit de montrer que 
$B$ est de dimension homologique finie comme $B\otimes_{K}B$-module (th\'eor\`eme de Serre), ou de mani\`ere \'equivalente de 
$Tor$-dimension finie. Or, par hypoth\`ese cette Tor-dimension est localement, pour la topologie $fppf$, finie sur 
$Spec\, B$, ce qui par descente $fppf$ de la platitude implique aussi qu'elle est localement finie pour la topologie
de Zariski. Par quasi-compacit\'e on conclut qu'elle est finie.
\hfill $\Box$ \\

\section{Le th\'eor\`eme de comparaison}

Le r\'esultat principal de ce travail est le th\'eor\`eme suivant.

\begin{thm}\label{t1}
Pour tout entier $n\geq 0$, le foncteur 
$$\phi_{n} : dSt_{et}^{n,li} \longrightarrow dSt_{fppf}^{n,pl}$$
est une \'equivalence de cat\'egories.
\end{thm}

Avant d'entrer dans les d\'etails de la preuve signalons quelques r\'eductions faciles.

\begin{lem}\label{l1}
Soit $n>0$.
\begin{enumerate}
\item Le foncteur $\phi_{n}$ est pleinement fid\`ele si tout $et$-champ d\'eriv\'e qui est $(n,li)$-g\'eom\'etrique est aussi un $fppf$-champ d\'eriv\'e. 

\item Si le foncteur $\phi_{n-1}$ est une \'equivalence alors le foncteur $\phi_{n}$ est pleinement fid\`ele.
\end{enumerate}
\end{lem}

\textit{Preuve:} $(1)$ Le foncteur $\phi_{n}$ est la restriction du foncteur $a$, de champ associ\'e pour la topologie $fppf$, qui poss\`ede comme adjoint \`a droite le foncteur d'inclusion $dSt_{fppf}(k) \subset dSt_{et}(k)$. Ainsi, $\phi_{n}$ est pleinement fid\`ele si pour tout $F\in dSt_{et}^{n,li}$, le morphisme d'adjonction 
$$F \longrightarrow a(F)$$
est un isomorphisme dans $dSt_{et}(k)$. Or, $F \longrightarrow a(F)$ est un isomorphisme pr\'ecis\'ement lorsque 
$F$ est un $fppf$-champ d\'eriv\'e. \\

$(2)$ Soit $F\in dSt_{et}^{n,li}$. D'apr\`es $(1)$ il nous suffit de montrer que le morphisme d'adjonction 
$f : F \longrightarrow \phi_{n}(F)$ est un isomorphisme dans $dSt_{et}(k)$. Soient $X$ et $Y$ deux objets affines et 
consid\'erons deux morphismes $X,Y \longrightarrow F$. Alors, en utilisant que $\phi_{n-1}$ est pleinement fid\`ele, on voit que 
le morphisme induit
$$X\times_{F}^{h}Y \longrightarrow \phi_{n-1}(X\times_{F}^{h}Y)\simeq X\times_{\phi_{n-1}(F)}^{h}Y$$
est un isomorphisme. Ceci implique que le morphisme $f$ est un monomorphisme (voir \cite[Rem. 1.2.6.2]{hagII}). Il nous reste donc \`a montrer que $f$ est aussi un \'epimorphisme de champs d\'eriv\'es pour la topologie \'etale. Soit $X$ un objet affine, et $x : X \longrightarrow \phi_{n}(F)$
un morphisme. On doit montrer qu'il existe un recouvrement \'etale $Y \longrightarrow X$ et un rel\`evement 
$Y \longrightarrow F$ du point $x$, c'est \`a dire que 
$$\xymatrix{
Y\ar[r] \ar[d] & F \ar[d]  \\
X \ar[r] & \phi_{n}(F)  }$$
commute dans $dSt_{et}$.
Soit $U=\coprod U_{i} \longrightarrow F$ un $n$-atlas lisse pour $F$. On consid\`ere le morphisme induit
$$\phi_{n}(U)\simeq U \longrightarrow \phi_{n}(F),$$
ainsi que 
$$U_{X}:=U\times^{h}_{\phi_{n}(F)}X \longrightarrow X.$$
On remarquera que $U_{X}$ est le produit fibr\'e homotopique de champs d\'eriv\'es \'etales qui sont aussi des champs d\'eriv\'es fppf, et donc est lui m\^eme un champs pour la topologie fppf.
Comme $\phi_{n}(F)$ est le $fppf$-champ associ\'e \`a $F$, il existe un recouvrement fppf $X' \longrightarrow X$ et un rel\`evement $X' \longrightarrow F$ de $x$.
En utilisant que $F \longrightarrow \phi_{n}(F)$ est un monomorphisme, on trouve
$$U_{X'}:=U_{X}\times^{h}_{X}X'\simeq U\times_{F}^{h}X',$$
ce qui montre que $U_{X'}$ est dans $dSt_{et}^{n-1,et}$. Soit $V=\coprod V_{i} \longrightarrow U_{X'}$
un $(n-1,li)$-atlas. Le compos\'e $V \longrightarrow U_{X}$ est un morphisme de $fppf$-champs d\'eriv\'es qui est $(n-2,pl)$-repr\'esentable, plat de pr\'esentation presque finie et surjectif. Ainsi, $U_{X}\in dSt_{fppf}^{n-1,pl}$, et comme $\phi_{n-1}$ est une \'equivalence on voit que 
$U_{X}\in dSt_{et}^{n-1,li}$. De plus, le morphisme $U \longrightarrow F$ \'etant lisse surjectif, le morphisme induit
$U_{X'} \longrightarrow X'$ est encore lisse surjectif. 
Comme le morphisme $X' \longrightarrow X$ est fid\`element plat cela implique que 
$U_{X} \longrightarrow X$ est lisse et surjectif. Ainsi, il existe un 
recouvrement \'etale $Y \longrightarrow X$ et un rel\`evement $Y \longrightarrow U_{X}$ de la projection 
$U_{X} \longrightarrow X$. En composant avec les morphismes
$U_{X} \longrightarrow U \longrightarrow F$, on trouve le rel\`evement du point $x$ cherch\'e.
\hfill $\Box$ \\

Ainsi, pour d\'emontrer le th\'eor\`eme \ref{t1} nous proc\'ederons par induction sur $n$. Pour $n=0$ l'\'enonc\'e est \'evident car
les champs d\'eriv\'es $0$-g\'eom\'etriques sont toujours les objets affines, et ce quelques soit le couple 
$(\tau,\mathbf{P})$ (toujours satisfaisant aux deux m\^emes conditions). 
On se fixe alors $n>0$, et on suppose que 
$\phi_{i}$ soit une \'equivalence pour tout $i<n$. Par le lemme \ref{l1} il nous suffit donc de d\'emontrer que 
$\phi_{n}$ est essentiellement surjectif. 
La d\'emonstration de cette derni\`ere assertion se fera en plusieurs \'etapes.

\subsection{Le champ d\'eriv\'e des extensions finies et strictement finies}

Dans cette section nous d\'efinissons des $et$-champs d\'eriv\'es $\underline{Fin}_{m}$ et 
$\underline{Fin}^{str}_{m}$, classifiant les sch\'emas d\'eriv\'es affines et de longueur finie fix\'ee $m$. L'objet 
$\underline{Fin}^{str}_{m}$ est une version rigidifi\'ee de $\underline{Fin}_{m}$, de sorte que $\underline{Fin}^{str}_{m}$ soit affine et qu'il existe un 
morphisme naturel $\underline{Fin}^{str}_{m} \longrightarrow \underline{Fin}_{m}$ qui soit un $Gl_{m}$-torseur. Cela implique en particulier que 
$\underline{Fin}_{m}$ est un $et$-champ d\'eriv\'e $(1,li)$-g\'eom\'etrique. \\

Pour commencer, nous dirons qu'un morphisme $A\longrightarrow B$ dans $sk-CAlg$ est 
\emph{fini et plat} si le $A$-module $B$ est projectif de type fini au sens de \cite[\S 1.2.4]{hagII}. Rappelons que cela signifie que 
$B$ est isomorphe, dans $Ho(sB-Mod)$ (la cat\'egorie homotopique des $B$-modules simpliciaux), \`a un r\'etracte d'un $A$-module de la forme $A^{m}$, pour un certain entier $m$.
En particulier, si $K$ est un corps et $A\longrightarrow K$ un morphisme, $B_{K}:=B\otimes_{A}^{\mathbb{L}}K$ est 
isomorphe, dans $Ho(sK-Mod)$ \`a un $K$-espace vectoriel de dimension finie. Nous dirons alors que 
$A \longrightarrow B$, plat fini, est de rang $m$ si pour tout corps $K$ et tout morphisme $A \longrightarrow K$ on a 
$dim_{K}(B\otimes_{A}^{\mathbb{L}}K)=m$. 

Par d\'efinition, $\underline{Fin}_{m}$ est le foncteur qui associe \`a $A\in sk-CAlg$ le nerf de la cat\'egorie des 
\'equivalences entre $A$-alg\`ebres plates, finies et de rang $m$. La construction pr\'ecise de ce foncteur
utilise, par exemple, la notion g\'en\'erale de pr\'efaisceaux de Quillen (voir \cite[App. B]{hagII}, ou aussi 
la fin du \S 2.3 de \cite{seat}). Concr\`etement le
foncteur
$$\underline{Fin}_{m} : sk-CAlg \longrightarrow SEns$$
se construit de la fa\c{c}on suivante. Pour $A\in sk-CAlg$ on consid\`ere $sA-CAlg^{fin,m,cof}$, la cat\'egorie dont les objets
sont les cofibrations de $k$-alg\`ebres simpliciales $A \longrightarrow B$, avec $B$ plat et fini, de rang $m$.
Les morphismes dans $sA-CAlg^{fin,m,cof}$ sont les diagrammes commutatifs 
$$\xymatrix{
B \ar[r]^-{f} & B' \\
 & A, \ar[ul] \ar[u]}$$
avec $f$ une \'equivalence faible. Si l'on a $A \longrightarrow A'$, un morphisme dans $sk-CAlg$, on dispose d'un 
foncteur de changement de bases
$$A'\otimes_{A} - : sA-CAlg^{fin,m,cof} \longrightarrow sA'-CAlg^{fin,m,cof}.$$
La construction $A \mapsto sA-CAlg^{fin,m,cof}$ d\'efinit de cette fa\c{c}on un pseudo-foncteur
$sk-CAlg \longrightarrow Cat$, que l'on rectifie, \`a \'equivalence pr\`es, par la construction de Grothendieck usuelle (voir par exemple \cite[App. B]{hagII}), afin d'obtenir un vrai foncteur 
$sk-CAlg \longrightarrow Cat$.
Compos\'e avec le foncteur nerf $Cat \longrightarrow SEns$, on obtient 
le foncteur cherch\'e
$$\underline{Fin}_{m} : sk-CAlg \longrightarrow SEns.$$
La valeur de ce foncteur sur un objet $A\in sk-CAlg$ est le nerf d'une cat\'egorie qui est naturellement \'equivalente
\`a $sA-CAlg^{fin,m,cof}$, et l'on dispose ainsi d'isomorphismes fonctoriels dans $Ho(SEns)$
$$\underline{Fin}_{m}(A)\simeq N(sA-CAlg^{fin,m,cof}).$$ 
On montre que le pr\'efaisceau simplicial 
$\underline{Fin}_{m}$ est un champ d\'eriv\'e pour la topologie fpqc, et en particulier pour les topologies 
\'etales et fppf (voir par exemple \cite[Thm. 1.3.7.2]{hagII} pour les grandes \'etapes de la preuve).

L'oubli de la structure multiplicative fournit, pour tout $A \in sk-CAlg$, un foncteur
$sA-CAlg^{fin,m,cof} \longrightarrow sA-Mod^{fin,m,cof}$, o\`u $sA-Mod^{fin,m,cof}$ est la cat\'egorie
des $A$-modules simpliciaux cofibrants, projectifs et de rang $m$. Cet oubli permet de construire un morphisme 
de $et$-champs d\'eriv\'es
$$\underline{Fin}_{m} \longrightarrow Vect_{m}\simeq BGl_{m},$$
o\`u $Vect_{m}$ est le champ des fibr\'es vectoriels de rang $n$ (voir \cite[Def. 1.3.7.5]{hagII}). La fibre homotopique de ce morphisme, prise au fibr\'e vectoriel trivial sera not\'ee $\underline{Fin}_{m}^{str}$.

\begin{prop}\label{p2}
Le $et$-champ d\'eriv\'e $\underline{Fin}_{m}^{str}$ est $0$-g\'eom\'etrique (i.e. affine), et le $et$-champ d\'eriv\'e
$\underline{Fin}_{m}$ est $(1,li)$-g\'eom\'etrique.
\end{prop}

\textit{Preuve:} Pour d\'emontrer cette proposition nous utiliserons le lemme de repr\'esentabilit\'e suivant
(voir \cite[App. C]{hagII} et \cite{lu} pour des versions plus g\'en\'erales).

\begin{lem}\label{l2}
Soit $F$ un $et$-champ d\'eriv\'e. On suppose que les conditions suivantes sont satisfaites.
\begin{enumerate}
\item Le $et$-champ tronqu\'e $t_{0}(F)$ est affine.

\item Le morphisme diagonal $F \longrightarrow F\times^{h} F$ est $0$-repr\'esentable. 

\item $F$ est nilcomplet: pour tout objet $A\in sk-CAlg$, de tour de Postnikov
$$A \longrightarrow \xymatrix{ \ar[r] & \dots A_{\leq k} \ar[r] & A_{\leq k-1} \ar[r] & 
\dots \ar[r] & A_{\leq 0}=\pi_{0}(A)}$$
le morphisme naturel 
$$F(A) \longrightarrow Holim_{k}F(A_{\leq k})$$
est une \'equivalence.

\item $F$ est inf-cart\'esien: pour tout $A\in sk-CAlg$, tout $A$-module connexe $M\in sA-Mod$ (i.e. $\pi_{0}(M)=0$), et toute
$k$-d\'erivation $d : A \longrightarrow A\oplus M$ (voir \cite[\S 1.2.1]{hagII}), le carr\'e homotopiquement cart\'esien suivant
$$\xymatrix{
A\oplus_{d}\Omega_{*}M \ar[r] \ar[d] & A \ar[d]^-{(id,0)} \\
A \ar[r]_-{d} & A\oplus M}$$
induit un carr\'e homotopiquement cart\'esien
$$\xymatrix{
F(A\oplus_{d}\Omega_{*}M) \ar[r] \ar[d] & F(A) \ar[d] \\
F(A) \ar[r] & F(A\oplus M).}$$

\end{enumerate}

Alors $F$ est affine.
\end{lem}

\textit{Preuve du lemme:} Tout d'abord, les conditions $(2)$ et $(4)$ impliquent que $F$ poss\`ede une th\'eorie de
l'obstruction au sens de \cite[\S 1.4.2]{hagII} (voir \cite[Prop. 1.4.2.7]{hagII}). On choisit $A_{0}$, une $k$-alg\`ebre commutative et un 
isomorphisme $\mathbb{R}\underline{Spec}\, A_{0}\longrightarrow t_{0}(F)$. En composant avec l'inclusion naturelle
$t_{0}(F) \longrightarrow F$ on trouve un morphisme
$$u_{0} : \mathbb{R}\underline{Spec}\, A_{0} \longrightarrow F.$$
Ce morphisme induisant un isomorphisme sur les tronqu\'es, on voit facilement  que son complexe
cotangent $\mathbb{L}_{u_{0}} \in Ho(sA_{0}-Mod)$ est $1$-connexe (i.e. $\pi_{0}(\mathbb{L}_{u_{0}})=
\pi_{1}(\mathbb{L}_{u_{0}})=0$). Nous allons construire par induction un diagramme commutatif
$$\xymatrix{
\mathbb{R}\underline{Spec}\, A_{0} \ar[d] \ar[r]^-{u_{0}} & F \\
\vdots \ar[d] & \\
\mathbb{R}\underline{Spec}\, A_{k} \ar[ruu]^-{u_{k}} \ar[d] & \\
\mathbb{R}\underline{Spec}\, A_{k+1} \ar[ruuu]_-{u_{k+1}} \ar[d] & \\
\vdots & }$$
o\`u $A_{k}$ est $k$-tronqu\'ee, le morphisme $A_{k+1} \longrightarrow A_{k}$ induit un 
isomorphisme sur les $\pi_{i}$ pour $i<k+1$, et de sorte \`a ce que
$\mathbb{L}_{u_{k}}$ soit $(k+1)$-connexe. On dispose donc d'un morphisme naturel de troncation
$\mathbb{L}_{u_{k}}\longrightarrow \pi_{k+2}(\mathbb{L}_{u_{k}})[k+2]$ (o\`u 
$(-)[m]$ d\'esigne le foncteur de suspension it\'er\'e $m$-fois).
Supposons que l'on ait construit $A_{k} \in sk-CAlg$, et un 
morphisme $u_{k} : \mathbb{R}\underline{Spec}\, A_{k} \longrightarrow F$ comme ci-dessus. On pose 
$M_{k+1}:=\pi_{k+2}(\mathbb{L}_{u_{k}})$, qui est un $\pi_{0}(A_{k})$-module. 
Le morphisme naturel $\mathbb{L}_{A_{k}} \longrightarrow \mathbb{L}_{u_{k}}$ induit 
donc un morphisme $\mathbb{L}_{A_{k}} \longrightarrow M_{k+1}[k+2]$, et donc une extension 
de carr\'e nul $A_{k+1}\longrightarrow A_{k}$, extension de $A_{k}$ par 
le $A_{k}$-module simplicial $\Omega_{*}M_{k}[k+2]\simeq M_{k}[k+1]$ (voir \cite[\S 1.2.1]{hagII}). Par construction de $A_{k+1}$, l'obstruction 
\`a \'etendre le morphisme $\mathbb{R}\underline{Spec}\, A_{k} \longrightarrow F$ le long de 
$\mathbb{R}\underline{Spec}\, A_{k+1} \longrightarrow \mathbb{R}\underline{Spec}\, A_{k}$ s'annule de mani\`ere naturelle  
et il existe donc une extension bien d\'efinie (voir \cite[Prop. 1.4.2.5]{hagII})
$$u_{k+1} : \mathbb{R}\underline{Spec}\, A_{k+1} \longrightarrow F.$$
Cette extension est telle que $\mathbb{L}_{u_{k+1}}$ est de plus $(k+2)$-connexe. Ceci finit de montrer l'existence
de la tour des morphismes $u_{k}$ comme ci-dessus. 

On pose alors $A:=holim_{k}A_{k}$. D'apr\`es la condition $(3)$, le syst\`emes des $u_{k}$ d\'efinit un 
morphisme $u : \mathbb{R}\underline{Spec}\, A \longrightarrow F$. Par construction ce morphisme induit un isomorphisme
sur les tronqu\'es $\mathbb{R}\underline{Spec}\, \pi_{0}(A) \simeq t_{0}(F)$, mais aussi un isomorphisme sur les 
complexes cotangents (i.e. $\mathbb{L}_{u}\simeq 0$). On d\'eduit de cela et de $(3)-(4)$, en utilisant les invariants de Postnikov (voir \cite[\S 2.2.1]{hagII}) que le morphisme $u$ induit, pour tout $B\in sk-CAlg$, une 
\'equivalence
$$(\mathbb{R}\underline{Spec}\, A)(B)=Map(A,B) \longrightarrow F(B).$$
On a donc bien $F\simeq \mathbb{R}\underline{Spec}\, A$. 
\hfill $\Box$ \\

Nous allons maintenant appliquer le lemme \ref{l2} au cas o\`u $F=\underline{Fin}_{m}^{str}$.
Soit $A$ une $k$-alg\`ebre commutative non simpliciale. On commence par remarquer 
que le foncteur $\pi_{0}$ induit un foncteur adjoint \`a gauche de la cat\'egorie
$sA-CAlg^{fin,m,cof}$ vers le groupo\"\i de des $A$-modules projectifs de type fini.
De cela on d\'eduit ais\'ement que 
le champ tronqu\'e $t_{0}(sA-CAlg^{fin,m,cof})$ est le champ usuel des
sch\'ema plats, finis et de longueur $m$. On d\'eduit de cela aussi que 
$t_{0}(\underline{Fin}^{str}_{m})$ n'est autre que le faisceaux des structures 
de $k$-alg\`ebres commutatives sur le module $k^{m}$, dont l'ensemble des valeurs sur une
$k$-alg\`ebre non simpliciale $A$ est l'ensemble des structures de $A$-alg\`ebres commuatives sur le
$A$-module $A^{m}$. Ce foncteur est clairement repr\'esentable par un sch\'ema affine. Pour pouvoir appliquer le lemme
\ref{l2} il nous faut donc d\'emontrer que les conditions $(2)$, $(3)$ et $(4)$ sont satisfaites pour 
$\underline{Fin}^{str}_{m}$. Comme $\underline{Fin}^{str}_{m}$ est la fibre homotopique du morphisme
$\underline{Fin}_{m} \longrightarrow Vect_{m}$, et que le champ d\'eriv\'e $Vect_{m}$ v\'erifie ces conditions (car il 
est $(1,li)$-g\'eom\'etrique), il nous suffit de montrer que 
$\underline{Fin}_{m}$ v\'erifie les conditions $(2)$, $(3)$ et $(4)$ du lemme \ref{l2}. \\

$(2)$ Montrer que le diagonale du $et$-champ d\'eriv\'e $\underline{Fin}_{m}$ est $0$-repr\'esentable est \'equivalent 
au fait suivant: soient $A \in sk-CAlg$, et $B$, $B'$ deux 
$A$-alg\`ebres simpliciales, cofibrantes et libres de rang $m$ comme $A$-modules. Alors le $et$-champ d\'eriv\'e
$\underline{Eq}(B,B')$ sur $\mathbb{R}Spec\, A$, qui envoie une $A$-alg\`ebre commutative $A'$ sur 
l'ensemble simplicial $\underline{Eq}(B\otimes_{A}A',B'\otimes_{A}A')$, des \'equivalences de 
$A'$-alg\`ebres commutatives, est affine. Pour cela on consid\`ere le $et$-champ d\'eriv\'e
$\underline{Hom}(B,B')$, de tous les morphismes entre $B$ et $B'$, et on commence par remarquer
qu'il est affine. En \'ecrivant $B$ comme une colimite homotopique de $A$-alg\`ebres simpliciales libres
on se ram\`ene au cas o\`u 
$B$ est une $A$-alg\`ebre commutative libre sur un ensemble $I$ 
(car les sch\'emas d\'eriv\'es affines sont stables par limites homotopiques). 
Dans ce cas
on a 
$$\underline{Hom}(B,B')\simeq \mathbb{R}\underline{Spec}\, Sym_{A}(\oplus_{I}(B')^{\vee}),$$
o\`u $(B')^{\vee}$ est le $A$-module simplicial dual de $B$. On remarque ensuite que 
$\underline{Eq}(B,B')$ est un ouvert Zariski de $\underline{Hom}(B,B')$. De plus, cet ouvert entre dans un carr\'e homotopiquement cart\'esien
$$\xymatrix{
\underline{Eq}(B,B') \ar[r] \ar[d] & \underline{Hom}(B,B') \ar[d] \\
Gl_{m} \ar[r] & M_{m},}$$
o\`u $M_{m}$ est le sch\'ema affine des matrices $m\times m$, et le morphisme
$\underline{Hom}(B,B') \longrightarrow M_{m}$ est obtenu en choisissant des isomorphismes dans $Ho(sA-Mod)$
de $B$ et $B'$ avec $A^{m}$. Ceci termine la preuve que la diagonale de $F$ est $0$-repr\'esentable. \\

$(3)$ En utilisant le point $(2)$, et le fait qu'un champ d\'eriv\'e affine v\'erifie les conditions
$(3)$ et $(4)$, on remarque que le morphisme naturel
$$F(A) \longrightarrow Holim_{k}F(A_{\leq k})$$
est un monomorphisme (i.e. est injectif sur les $\pi_{0}$ et un isomorphisme sur tous les
$\pi_{i}$ pour $i>0$). Il reste donc \`a voir que ce morphisme est surjectif sur les 
composantes connexes. Il est facile de voir qu'il existe une bijection 
entre $\pi_{0}(Holim_{k}F(A_{\leq k}))$, et l'ensemble des classes d'isomorphismes 
d'objets dans $Holim_{k}Ho(sA_{\leq k}-CAlg^{fin,m,cof})$. Ainsi, un \'el\'ement
de $\pi_{0}(Holim_{k}F(A_{\leq k}))$ se repr\'esente par un syst\`eme d'objets
$B_{k}\in sA_{\leq k}-CAlg^{fin,m,cof}$, et des \'equivalences 
$B_{k+1}\otimes_{A_{\leq k+1}}A_{\leq k} \longrightarrow B_{k}$, de $A_{\leq k}$-alg\`ebres
simpliciales commutatives. A un tel syst\`eme on associe 
$$B:=Holim_{k}B_{k} \in Ho(sA-CAlg).$$
Il n'est pas difficile de voir que $B$ est fini de rang $m$ comme $B$-module simplicial, et qu'un remplacement 
cofibrant de $B$ fournit un ant\'ec\'edent, \`a homotopie pr\`es, du syst\`eme $\{B_{k}\}$
par le morphisme $F(A) \longrightarrow Holim_{k}F(A_{\leq k})$. \\

$(4)$ Tout comme pour le point $(3)$, le morphisme en question est un monomorphisme. Pour montrer qu'il est aussi surjectif sur les composantes connexes on utilise le lemme \cite[Lem. 4.2]{dhall} (ou plus pr\'ecis\`ement sa version pour des sous-cat\'egories de mod\`eles stables par \'equivalence). Nous laissons le lecteur v\'erifier que le
carr\'e de foncteur de Quillen
$$\xymatrix{
s(A\oplus_{d}\Omega_{*}M)-CAlg \ar[r] \ar[d] & sA-CAlg \ar[d]^-{(id,0)} \\
sA-CAlg \ar[r]_-{d} & s(A\oplus M)-CAlg
}$$
induit un carr\'e homotopiquement cart\'esien
$$\xymatrix{
N(s(A\oplus_{d}\Omega_{*}M)-CAlg^{fin,m,cof}) \ar[r] \ar[d] & N(sA-CAlg^{fin,m,cof}) \ar[d]^-{(id,0)} \\
N(sA-CAlg^{fin,m,cof}) \ar[r]_-{d} & N(s(A\oplus M)-CAlg^{fin,m,cof}).
}$$

Ceci termine la v\'erification que le champ d\'eriv\'e $Fin^{m}$ v\'erifie les
conditions du lemme \ref{l2}, et donc la preuve de la proposition \ref{p2}.
\hfill $\Box$ \\

Le complexe cotangent du champ $\underline{Fin}_{m}$ peut se d\'ecrire de la fa\c{c}on suivante. Il s'agit du fait, bien connu, que les d\'eformations infinit\'esimales d'une $A$-alg\`ebre commutative projective et de rang finie 
sont d\'ecrites par la cohomologie d'Andr\'e-Quillen.

\begin{prop}\label{p3}
Soit $A\in sk-CAlg$, et 
$u : \mathbb{R}\underline{Spec}\, A \longrightarrow \underline{Fin}_{m}$ un morphisme de champs d\'eriv\'es correspondant \`a 
une $A$-alg\`ebre commutative $B$ projective et de rang $m$ comme $A$-module. Notons
$$B^{\vee}:=\mathbb{R}\underline{Hom}_{sA-Mod}(B,A)$$ 
le $A$-module dual de $B$, muni de sa structure de $B$-module naturelle.
Alors,  il existe un isomorphisme naturel dans $Ho(Sp(sA-Mod))$
$$\mathbb{L}_{\underline{Fin}_{m},u} \simeq \mathbb{L}_{B/A}[-1]\otimes_{B}^{\mathbb{L}}B^{\vee}.$$
\end{prop}

\textit{Preuve:} Soit $\Omega_{u}\underline{Fin}_{m}$, le champ d\'eriv\'e des lacets de base $u$ dans $\underline{Fin}_{m}$. Comme nous l'avons d\'ej\`a vu lors de la preuve de la proposition \ref{p1}, le champ d\'eriv\'e $\Omega_{u}\underline{Fin}_{m}$ est un ouvert
du champ $\underline{Hom}(B,B)$, des endomorphismes de $B$ comme $A$-alg\`ebre simpliciale. On a ainsi
$$\mathbb{L}_{\underline{Fin}_{m}}[1] \simeq \mathbb{L}_{\Omega_{u}\underline{Fin}_{m},u} \simeq \mathbb{L}_{\underline{Hom}(B,B),id}.$$
Il est facile de voir, par d\'efinition du complexe cotangent et du champ 
d\'eriv\'e $\underline{Hom}(B,B)$, qu'il existe des isomorphismes fonctoriel en $M\in Ho(sA-Mod)$ 
$$[\mathbb{L}_{\underline{Hom}(B,B),id},M]_{Ho(sA-Mod)}\simeq [\mathbb{L}_{B/A},M\otimes_{A}^{\mathbb{L}}B]_{Ho(sB-Mod)}$$
$$\simeq [\mathbb{L}_{B/A},\mathbb{R}\underline{Hom}_{sA-Mod}(B^{\vee},M)]_{Ho(sB-Mod)}
\simeq [\mathbb{L}_{B/A}\otimes_{B}^{\mathbb{L}}B^{\vee},M]_{Ho(sA-Mod)},$$
ce qui implique l'\'enonc\'e de la proposition.
\hfill $\Box$ \\

\subsection{Restriction des scalaires le long d'un morphisme fini}

Soit $p : G \longrightarrow H$ un morphisme de $et$-champs d\'eriv\'es, et consid\'erons les 
cat\'egories de champs au-dessus de $G$ et de $H$
$$dSt_{et}(k/F):=Ho(dAff_{k}^{\sim,et}/G) \qquad 
dSt_{et}(k/F'):=Ho(dAff_{k}^{\sim,et}/H).$$
Le changement de bases le long du morphisme $p$ induit un foncteur
$$p^{*} :  dSt_{et}(k/H)  \longrightarrow  dSt_{et}(k/G)$$
qui envoie un objet $(F \rightarrow H)$ sur $(F\times_{H}^{h}G \rightarrow G)$.
Ce foncteur poss\`ede un adjoint \`a droite
$$p_{*} : dSt_{et}(k/G) \longrightarrow dSt_{et}(k/H),$$
qui \`a un objet $F \rightarrow G$ associe le $et$-champ d\'eriv\'e
des morphismes au-dessus de $H$
$$p_{*}(F):=\mathbf{Map}_{/H}(G,F),$$
o\`u $\mathbf{Map}_{/H}$ d\'esigne le Hom interne de la cat\'egorie
$dSt_{et}(k/H)$ (voir \cite[\S 3.6]{hagI} pour l'existence des Hom internes). Le foncteur
$p_{*}$ s'appelle la \emph{restriction des scalaires le long du morphisme $p$}. 

L'adjonction $(p^{*},p_{*})$ peut aussi se r\'ealiser 
par une adjonction de Quillen de la fa\c{c}on suivante. On repr\'esente $p$, \`a \'equivalence pr\`es, par une fibration $p : G \longrightarrow H$ entre objets fibrants de $dAff_{k}^{\sim,et}$. On remarque alors que le foncteur de changement de base 
$$-\times_{H}G : dAff_{k}^{\sim,et}/H \longrightarrow dAff_{k}^{\sim,et}/G$$
est de Quillen \`a gauche. L'adjonction induite sur les cat\'egories homotopiques correspondantes est l'adjonction d\'ecrite ci-dessus
$(p^{*},p_{*})$. \\

L'\'enonc\'e suivant est un cas particulier d'un crit\`ere de repr\'esentabilit\'e pour la
restriction des scalaires le long de morphismes propres et plats (d\'emontr\'e par exemple dans \cite{lu}) 
(se restreindre aux morphismes 
finis simplifie sensiblement la preuve).

\begin{prop}\label{p4}
Soit $p : G \longrightarrow H$ un morphisme entre $et$-champs d\'eriv\'es. On suppose que 
$p$ est $0$-g\'eom\'etrique, plat et fini: pour tout $X=\mathbb{R}\underline{Spec}\, A \longrightarrow H$, 
on a $G\times_{H}^{h}X \simeq \mathbb{R}\underline{Spec}\, B$ avec $B$ une $A$-alg\`ebre commutative projective et finie. Alors le foncteur
$$p_{*} : dSt_{et}(k/G) \longrightarrow dSt_{et}(k/H)$$
pr\'es\`erve les champs d\'eriv\'es $(n,li)$-repr\'esentables, pour tout $n\geq 0$.
\end{prop}

\textit{Preuve:} Tout d'abord, \cite[Prop. 1.3.3.4]{hagII} implique que notre assertion 
est locale pour la topologie \'etale sur $H$ 
et l'on peut donc supposer que $G$ et $H$ sont tous deux affines, et que le morphisme 
$$p : G=\mathbb{R}\underline{Spec}\, B \longrightarrow H=\mathbb{R}\underline{Spec}\, A$$
correspond \`a une $A$-alg\`ebre commutative $B$, libre et de rang fini comme $A$-module. Les cat\'egories
$dSt_{et}(k/G)$ et $dSt_{et}(k/H)$ sont alors respectivement \'equivalentes \`a 
$dSt_{et}(B)$ et $dSt_{et}(A)$, des champs pour la topologie \'etales sur 
$dAff_{B}:=(sB-CAlg)^{op}$ et $dAff_{A}:=(sA-CAlg)^{op}$. 
Le foncteur $p_{*}$ est alors donn\'e par la formule
$$p_{*}(F)(A'):=F(A'\otimes_{A}^{\mathbb{L}}B) \qquad \forall \, A'\in sA-CAlg.$$

La preuve proc\`ede, comme il se doit, par induction sur $n$. Commen\c{c}ons par traiter le cas $n=0$.
Soit $F=\mathbb{R}\underline{Spec}\, B'$, affine au-dessus de $\mathbb{R}\underline{Spec}\, B$. On a
$$p_ {*}(F)(A'):=F(A'\otimes_{A}^{\mathbb{L}}B)\simeq \mathbb{R}\underline{Hom}_{sB-CAlg}(B',A'\otimes_{A}^{\mathbb{L}}B).$$
Pour montrer que $p_{*}(F)$ est affine, nous utiliserons le r\'esultat de repr\'esentabilit\'e des foncteurs d\'efinis sur des cat\'egories de mod\`eles combinatoires \cite[Prop. 1.9]{tv}. Il nous suffit donc de montrer les
deux assertions suivantes.
\begin{enumerate}
\item  Le foncteur $p_{*}(F)$ commute avec les limites homotopiques dans $sA-CAlg$.
\item Le foncteur $p_{*}(F)$ commute avec les colimites $\lambda$-filtrantes pour 
$\lambda$ un cardinal suffisemment grand. 
\end{enumerate}

Le point $(1)$ se d\'eduit du fait que $B$ est libre de rang fini comme $A$-module. En effet, il suffit de voir que $A' \mapsto A'\otimes_{A}^{\mathbb{L}}B$ commute aux limites homotopiques. Comme le foncteur 
d'oubli $sB-CAlg \longrightarrow sB-Mod$ refl\`ete les limites homotopiques il suffit de voir que 
le foncteur 
$$-\otimes_{A}^{\mathbb{L}}B : Ho(sA-Mod) \longrightarrow Ho(sB-Mod)$$
commute aux limites homotopiques. Mais, comme $B\simeq A^{m}$ comme $A$-module, ce dernier foncteur
est isomorphe \`a $M \mapsto M^{m}$, et commute donc aux limites homotopiques. 

Le point $(2)$ se d\'eduit du fait que la cat\'egorie de mod\`eles $sB-CAlg$ est combinatoire. Dans une telle cat\'egorie
de mod\`eles tout objet $x$ est homotopiquement $\lambda$-petit pour un cardinal $\lambda$ (i.e.
$Map(x,-)$ commute aux colimites homotopiques $\lambda$-filtrantes, voir \cite{du}). \\

Cela termine la preuve de la proposition \ref{p3} pour le cas $n=0$. Supposons maintenant qu'elle soit d\'emontr\'ee pour $n$ et montrons qu'elle reste vraie au rang $n+1$. Soit $F\in dSt_{et}(A)$ un 
$et$-champ d\'eriv\'e $(n+1,li)$-g\'eom\'etrique. Soit $\{U_{i}\}$ des affines et 
$$f : U:=\coprod_{i}U_{i} \longrightarrow F$$
un $(n+1,li)$-atlas. Nous allons montrer, par r\'ecurrence sur $n$, que le morphisme induit
$$g : V:=\coprod_{i}p_{*}(U_{i}) \longrightarrow p_{*}(F)$$
est un $(n+1,li)$-atlas. Nous supposerons donc que $p_{*}$ pr\'eserve les $(m,li)$-atlas pour
$m\geq n$. 

Nous avons d\'ej\`a vu que chacun des $p_{*}(U_{i})$ \'etait affine (cas $n=0$). Commen\c{c}ons par voir que le morphisme
$g$ est un \'epimorphisme de champs d\'eriv\'es pour la topologie \'etale, et plus g\'en\'eralement que le foncteur $p_{*}$ pr\'eserve les \'epimorphismes pour la topologie \'etale. Pour voir cela, soit 
$G \longrightarrow G'$ un \'epimorphisme dans $dSt_{et}(B)$, et soit $x\in p_{*}(G')(A')$ pour 
$A'\in sA-CAlg$. Comme $p_{*}(G')(A')\simeq G'(A'\otimes_{A}^{\mathbb{L}}B)$, il existe un 
recouvrement \'etale $A'\otimes_{A}^{\mathbb{L}}B \longrightarrow C$ tel que $x$ se rel\`eve, \`a homotopie pr\`es, \`a $G(C)$. Or, comme $B$ est finie et plat sur $A$, on sait qu'il existe un recouvrement \'etale
$A' \longrightarrow C'$, et un diagramme commutatif dans $Ho(sk-CAlg)$
$$\xymatrix{
A'\otimes_{A}^{\mathbb{L}}B \ar[r] \ar[d] & C \ar[dl] \\
C'\otimes_{A}^{\mathbb{L}}B. & }$$
(pour d\'emontrer l'existence d'un tel $C'$ on utilise \cite[Cor. 2.2.2.9]{hagII} 
pour ramener l'\'enonc\'e au cas 
des $k$-alg\`ebres commutatives non simpliciales, et on utilise qu'une alg\`ebre finie sur 
un anneau local henselien strict et un produit d'anneaux locaux hens\'eliens stricts).
Cela montre que $x$ se rel\`eve, \`a homotopie pr\`es, \`a un \'el\'ement dans 
$G(C'\otimes_{A}^{\mathbb{L}}B)\simeq p_{*}(G)(C')$, ce qu'il fallait montrer. 

Comme le morphisme $g$ est un \'epimorphisme, le champ d\'eriv\'e $p_{*}(F)$ est \'equivalent au quotient 
homotopique du groupoide de Segal nerf de $g$ 
$$p_{*}(F) \simeq Hocolim \left( [n] \mapsto \underbrace{V\times_{p_{*}(F)}^{h}V\times_{p_{*}(F)}V \dots {\times}_{p_{*}(F)}^{h}V}_{n-fois} \right)$$
(voir \cite[\S 1.3.4]{hagII}). Ainsi, pour voir que 
$p_{*}(F)$ est $(n+1,li)$-g\'eom\'etrique il suffit de voir que ce nerf est un groupoide de Segal 
$(n,li)$-g\'eom\'etrique et lisse, c'est \`a dire la projection 
$$V\times_{p_{*}(F)}^{h}V \longrightarrow V$$
est un morphisme $(n,li)$-repr\'esentable et lisse.
Comme le morphisme $p_{*}$ commute aux limites homotopiques (car c'est un foncteur d\'eriv\'e \`a droite d'un foncteur de Quillen \`a droite) on a 
$$V\times_{p_{*}(F)}^{h}V \simeq \coprod_{i,j}p_{*}(U_{i}\times^{h}_{F}U_{j}),$$
et par induction on voit que la projection $V\times_{p_{*}(F)}^{h}V \longrightarrow V$ est $(n,li)$-repr\'esentable. Pour terminer la preuve de la proposition il nous reste \`a montrer que cette projection est lisse, c'est \`a dire que pour tout $i,j$ la projection
$$p_{*}(U_{i}\times^{h}_{F}U_{j}) \longrightarrow p_{*}(U_{i})$$
est lisse. Pour cela nous utiliserons le crit\`ere infinit\'esimal de lissit\'e 
\cite[Prop. 2.2.5.1]{hagII}. Tout d'abord, 
le fait que $p_{*}(U_{i}\times^{h}_{F}U_{j}) \longrightarrow p_{*}(U_{i})$ soit localement de pr\'esentation finie se d\'eduit par induction sur $n$ et par le fait que $p_{*}$ pr\'eserve les objets affines
et de pr\'esentation finie. 

Soit donc $A'\in sA-CAlg$, $M$ un $A'$-module connexe et $d : A' \longrightarrow A'\oplus M$ une 
$A$-d\'erivation. Notons, $B':=A'\otimes_{A}^{\mathbb{L}}B$, $M_{B}:=M\otimes_{A}^{\mathbb{L}}B$, ainsi que 
$$d_{B} : B'  \longrightarrow B'\oplus M_{B}$$
la $B$-d\'erivation induite par changement de base le long de $A\rightarrow B$.
Soit $A'\oplus_{d}\Omega_{*}M$ l'extension de carr\'e nul associ\'ee \`a $d$, et consid\'erons le morphisme
$$p_{*}(U_{i}\times^{h}_{F}U_{j})(A'\oplus_{d}\Omega_{*}M) \longrightarrow
p_{*}(U_{i})(A'\oplus_{d}\Omega_{*}M)\times_{p_{*}(U_{i})(A')}^{h}p_{*}(U_{i}\times^{h}_{F}U_{j})(A').$$
Il s'agit de montrer que ce morphisme est surjectif \`a homotopie pr\`es. Tout d'abord, comme $B$ est plat sur $A$ on a 
$$(A'\oplus_{d}\Omega_{*}M)\otimes_{A}^{\mathbb{L}}B \simeq
B'\oplus_{d_{B}}\Omega_{*}(M_{B}).$$
Ainsi, 
par adjonction le morphisme ci-dessus 
s'\'ecrit aussi
$$(U_{i}\times^{h}_{F}U_{j})(B'\oplus_{d_{B}}\Omega_{*}M_{B}) \longrightarrow
U_{i}(B'\oplus_{d_{B}}\Omega_{*}M_{B})\times_{U_{i}(B')}^{h} (U_{i}\times^{h}_{F}U_{j})(B').$$
Or ce morphisme est surjectif sur les composantes connexes car 
$U_{i}\times^{h}_{F}U_{j} \longrightarrow U_{i}$ est un morphisme $(n,li)$-repr\'esentable et lisse
(voir \cite[Prop. 2.2.5.1]{hagII}). 
Ceci termine la preuve de la proposition \ref{p3}. 
\hfill $\Box$ \\

\subsection{Quasi-sections}

Dans ce paragraphe nous introduisons le $et$-champ d\'eriv\'e $\underline{QSect_{f}}$, des \emph{quasi-sections} d'un morphisme fix\'e $f : F' \longrightarrow F$ dans $dSt_{et}(k)$. Les quasi-sections sont des sections \`a extension finie pr\`es: 
le champ $\underline{QSect_{f}}$ classifie les diagrammes commutatifs
$$\xymatrix{
X' \ar[r] \ar[d] & F' \ar[d]_-{f} \\
X \ar[r] & F}$$
avec $X' \longrightarrow X$ un morphisme $0$-repr\'esentable, plat et fini. La d\'efinition pr\'ecise du champ 
d\'eriv\'e $\underline{QSect_{f}}$ va occuper toute la premi\`ere partie de ce paragraphe. Nous \'etudierons ensuite la repr\'esentabilit\'e de $\underline{QSect_{f}}$ en fonction de celle de $p$, ce qui est une fa\c{c}on de combiner les r\'esultats \ref{p2} et \ref{p3}. Nous introduirons aussi un certain sous-champ ouvert des quasi-sections \emph{quasi-lisses}. \\

Nous commencerons par le $et$-champ d\'eriv\'e des carr\'es commutatifs, not\'e $\underline{Car}$, et 
d\'efini de la fa\c{c}on suivante. On note $\square$ la cat\'egorie classifiant les carr\'es commutatifs
$$\square:=\Delta^{1} \times \Delta^{1},$$
o\`u $\Delta^{1}$ est la cat\'egorie avec deux objets $0$ et $1$, et un unique morphisme 
$0 \rightarrow 1$
$$\Delta^{1}=\left( 0 \rightarrow 1\right) \qquad
\square=\left( 0 \rightarrow 1\right) \times \left( 0 \rightarrow 1\right).$$
Pour $A \in sk-CAlg$, on consid\`ere la cat\'egorie
de mod\`eles $dAff_{A}^{\sim,et}$, ainsi que la cat\'egorie des diagrammes
$(dAff_{A}^{\sim,et})^{\square}$. On d\'efinit $\underline{Car}(A)$ comme \'etant le nerf de la cat\'egorie des \'equivalences dans $(dAff_{A}^{\sim,et})^{\square}$
$$\underline{Car}(A):=N(w(dAff_{A}^{\sim,et})^{\square}).$$
Lorsque $A \rightarrow B$ est un morphisme dans $sk-CAlg$, on dispose d'un foncteur d'oubli
$dAff_{B} \longrightarrow dAff_{A}$, qui induit un foncteur de restriction
$(dAff_{A}^{\sim,et})^{\square} \longrightarrow (dAff_{B}^{\sim,et})^{\square}$. Ce foncteur de restriction
pr\'eserve les \'equivalences et induit donc un morphisme d'ensembles simpliciaux
$\underline{Car}(A) \longrightarrow \underline{Car}(B)$. L'association $A \mapsto 
\underline{Car}(A)$ d\'efinit de cette fa\c{c}on un pr\'efaisceau simplicial $\underline{Car}$ sur
$dAff_{k}$. On v\'erifie que $\underline{Car}$ est un $et$-champ d\'eriv\'e (en utilisant par exemple
les techniques pr\'esent\'ees \`a la fin du \S 2.3 de \cite{seat}).

De la m\^eme fa\c{c}on, on d\'efinit $\underline{Mor}$, le $et$-champ d\'eriv\'e des morphismes, qui envoie
$A$ sur le nerf des \'equivalences dans $(dAff_{A}^{\sim,et})^{\Delta^{1}}$. \\

Fixons maintenant $f : F' \longrightarrow F$ un morphisme dans $dSt_{et}(k)$, que nous prenons soins de relever en un morphisme dans $dAff_{k}^{\sim,et}$ (quitte \`a prendre des remplacements fibrants et cofibrants). L'inclusion en l'objet $1$
$$\{1\}\times id : \Delta^{1} \hookrightarrow \square=\Delta^{1}\times \Delta^{1}$$
induit un morphisme de restriction
$$\underline{Car} \longrightarrow \underline{Mor}.$$
Ce morphisme envoie un carr\'e commutatif
$$\xymatrix{
G' \ar[r] \ar[d] & F' \ar[d] \\
G \ar[r] & F}$$
sur le morphisme $F' \rightarrow F$. Le morphisme $f$ d\'efinit un point global 
$* \longrightarrow \underline{Mor}$, et on pose
$$\underline{Car}_{/f}:=*\times^{h}_{\underline{Mor}}\underline{Car}.$$
Le $et$-champ d\'eriv\'e $\underline{Car}$ est le champ des carr\'es commutatifs dont le morphisme but est fix\'e \'egal \`a $f$. 

On consid\`ere un sous-champ $\underline{Car}_{m}' \subset \underline{Car}$: pour $A \in sk-CAlg$, 
$\underline{Car}_{m}'(A)$ est la r\'eunion des composantes connexes de $\underline{Car}(A)$
form\'ee des diagrammes commutatifs dans $dSt_{et}(A)$
$$\xymatrix{
G' \ar[r] \ar[d] & F' \ar[d] \\
G \ar[r] & F}$$
avec $G\simeq \mathbb{R}\underline{Spec}\, A$ et $G' \longrightarrow \mathbb{R}\underline{Spec}\, A$ 
plat et fini de rang $m$ (c'est un sous-champ car 
\^etre affine est une condition locale, ainsi qu'\^etre plat et fini). 
On consid\`ere alors le
carr\'e homotopiquement cart\'esien suivant
$$\xymatrix{
\underline{QSect}_{f,m} \ar[r] \ar[d] & \underline{Car}_{/f} \ar[d] \\
\underline{Car}_{m}' \ar[r] & \underline{Car}.}$$

\begin{df}\label{d1}
Soit $f : F' \longrightarrow F$ comme ci-dessus.
Le $et$-\emph{champ d\'eriv\'e des quasi-sections de $f$ de rang $m$} est 
$\underline{QSect}_{f,m} \in dSt_{et}(k)$ d\'efini ci-dessus.
\end{df}

Notons qu'il existe une projection naturelle $\underline{QSect}_{f,m} \longrightarrow F$, qui 
\`a un diagramme
$$\xymatrix{
G' \ar[r] \ar[d] & F' \ar[d] \\
G \ar[r] & F}$$
avec $G\simeq \mathbb{R}\underline{Spec}\, A$, associe le $A$-point $x$ de $F$ correspondant
(nous laissons au lecteur le soin de d\'efinir cette projection de mani\`ere rigoureuse). \\

On combine maintenant les deux r\'esultats de repr\'esentabilit\'e Prop. \ref{p2} et \ref{p4}
en l'\'enonc\'e suivant.

\begin{prop}\label{p5}
Si $f : F' \longrightarrow F$ est un morphisme $(n,li)$-repr\'esentable dans $dSt_{et}(k)$, avec $n>0$. Alors la projection
$\underline{QSect}_{f,m} \longrightarrow F$ est $(n,li)$-g\'eom\'etrique.
\end{prop}

\textit{Preuve:} On dispose d'un morphisme $\underline{QSect}_{f,m} \longrightarrow \underline{Fin}_{m}$, 
qui \`a un diagramme commutatif
$$\xymatrix{
G' \ar[r] \ar[d] & F' \ar[d] \\
\mathbb{R}\underline{Spec}\, A \ar[r] & F}$$
associe le champ affine $G'$, de la forme $\mathbb{R}\underline{Spec}\, B$, pour $B$ une $A$-alg\`ebre
plate et finie de rang $m$ (nous laissons le soin au lecteur de revenir sur les d\'efinition 
de $\underline{QSect}_{f,m}$ et $\underline{Fin}_{m}$ afin de d\'efinir proprement ce morphisme
$\underline{QSect}_{f,m} \rightarrow \underline{Fin}_{m}$ dans $dSt_{et}(k)$). Si l'on se fixe
$A \rightarrow B$ plate et finie de rang $m$, on dispose d'un carr\'e homotopiquement cart\'esien
$$\xymatrix{
\underline{Map}_{/\mathbb{R}\underline{Spec}\, A}(\mathbb{R}\underline{Spec}\, B,F'\times_{F}^{\mathbb{L}}\mathbb{R}\underline{Spec}\, A) \ar[r] \ar[d] & 
\underline{QSect}_{f,m} \ar[d] \\
\mathbb{R}\underline{Spec}\, A \ar[r]_-{B} & \underline{Fin}_{m}\times^{h}F.}$$
La proposition \ref{p4}, et l'hypoth\`ese de repr\'esentabilit\'e sur $f$, implique que
$\underline{Map}_{/\mathbb{R}\underline{Spec}\, A}(\mathbb{R}\underline{Spec}\, B,F'\times_{F}^{\mathbb{L}}\mathbb{R}\underline{Spec}\, A)$ est 
$(n,li)$-repr\'esentable au-dessus de $\mathbb{R}\underline{Spec}\, A$. Ceci montre que le morphisme 
$\underline{QSect}_{f,m} \rightarrow \underline{Fin}_{m}\times^{h}F$ est $(n,li)$-g\'eom\'etrique. Comme
$\underline{Fin}_{m}$ est $(1,li)$-g\'eom\'etrique (prop. \ref{p2}) il s'en suit 
que $\underline{QSect}_{f,m} \longrightarrow F$ est $(n,li)$-repr\'esentable (on utilise ici $n>0$).
\hfill $\Box$ \\

Nous aurons besoin d'un r\'esultat un peu plus fin. Consid\'erons le carr\'e homotopiquement cart\'esien suivant
$$\xymatrix{
\underline{QSect}_{f,m}^{str} \ar[r] \ar[d] & \underline{QSect}_{f,m} \ar[d] \\
\underline{Fin}_{m}^{str} \ar[r] & \underline{Fin}_{m}. }$$

\begin{df}\label{d2}
Soit $f : F' \longrightarrow F$ comme ci-dessus.
Le $et$-\emph{champ d\'eriv\'e des quasi-sections strictes de $f$ de rang $m$} est 
$\underline{QSect}^{str}_{f,m} \in dSt_{et}(k)$ d\'efini ci-dessus.
\end{df}

\begin{prop}\label{p6}
Si $f : F' \longrightarrow F$ est un morphisme $(n,li)$-repr\'esentable dans $dSt_{et}(k)$. Alors
la projection naturelle $\underline{QSect}_{f,m}^{str} \longrightarrow F$ est $(n,li)$-g\'eom\'etrique.
\end{prop}

\textit{Preuve:} Au cours de la preuve de la proposition \ref{p5} nous avons vu que le morphisme
$\underline{QSect}_{f,m} \longrightarrow \underline{Fin}_{m}\times^{h}F$ \'etait $(n,li)$-repr\'esentable (ce point n'utilisait pas l'hypoth\`ese $n>0$). Il s'en suit par changement de bases que 
$\underline{QSect}_{f,m}^{str} \longrightarrow \underline{Fin}_{m}^{str}\times^{h}F$ est $(n,li)$-repr\'esentable. 
Comme $\underline{Fin}_{m}^{str}$ est affine, on trouve que $\underline{QSect}_{f,m}^{str} \longrightarrow F$ est 
$(n,li)$-g\'eom\'etrique. \hfill $\Box$ \\

Pour terminer nous allons consid\'erer un certain sous-champ
$\underline{QSect}_{f,m}^{str,ql} \subset \underline{QSect}_{f,m}^{str}$, form\'ee des quasi-sections
\emph{quasi-lisses}. Nous continuons avec un morphisme $f : F' \longrightarrow F$ que nous supposons
$(n,li)$-repr\'esentable. Soit $A \in sk-CAlg$, et 
soit $q : \mathbb{R}\underline{Spec}\, A \longrightarrow \underline{QSect}_{f,m}^{str}$ un morphisme
correspondant \`a un carr\'e commutatif
$$\xymatrix{
\mathbb{R}\underline{Spec}\, B \ar[r]^-{u} \ar[d] & F' \ar[d] \\
\mathbb{R}\underline{Spec}\, A \ar[r] & F}$$
Nous dirons que la quasi-section $q$ est \emph{quasi-lisse} (on pourrait aussi dire \emph{l.c.i.}) si le complexe cotangent
$\mathbb{L}_{u}$ du 
morphisme $u$ est parfait et d'amplitude contenue dans $[-1,\infty[$ (on rappelle que $\mathbb{L}_{u}$ est la cofibre homotopique du morphisme $\mathbb{L}_{F,u} \rightarrow \mathbb{L}_{B}$). 
Nous noterons $\underline{QSect}_{f,m}^{str,ql}$ le sous-pr\'efaisceau simplicial de
$\underline{QSect}_{f,m}^{str}$ form\'e des quasi-sections qui sont quasi-lisses: pour 
$A\in sk-CAlg$, l'ensemble simplicial $\underline{QSect}_{f,m}^{str,ql}(A)$ est 
la r\'eunion des composantes connexes de $\underline{QSect}_{f,m}^{str}(A)$ qui consistent 
est des carr\'es comme ci-dessus avec $u$ quasi-lisse.

\begin{prop}\label{p7}
Soit $f : F' \longrightarrow F$ un morphisme $(n,li)$-repr\'esentable plat et presque de pr\'esentation finie.
Le morphisme d'inclusion 
$$\underline{QSect}_{f,m}^{str,ql} \longrightarrow \underline{QSect}_{f,m}^{str}$$
est une immerison de Zariski ouverte. 
\end{prop}

\textit{Preuve:} Un morphisme $A \longrightarrow B$ dans $sk-CAlg$ est dit \emph{quasi-lisse} s'il 
est homotopiquement de pr\'esentation finie et si de plus son complexe cotangent
$\mathbb{L}_{B/A}$ est homotopiquement de pr\'esentation finie dans $sB-Mod$ (on dit aussi \emph{parfait}, tout au moins lorsque l'on consid\`ere $\mathbb{L}_{B/A}$ comme un $B$-module \emph{stable}, 
voir \cite[\S 1.2.11]{hagII}), 
et d'amplitude contenue dans $[-1,0]$. Reppelons que cette derni\`ere condition signifie que pour tout $B$-module connexe et simplement connexe
$M$ on a $[\mathbb{L}_{B/A},M]=0$. En utilisant \cite[Prop. 2.2.2.4]{hagII} on voit qu'un tel morphisme est 
quasi-lisse si et seulement si $\pi_{0}(B)$ est une $\pi_{0}(A)$-alg\`ebre de pr\'esentation finie et si de plus
$\mathbb{L}_{B/A}$ est parfait et d'amplitude contenue dans $[-1,0]$. On remarque que la notion 
de quasi-lissit\'e est locale pour la topologie \'etale sur $dAff_{k}$, et s'\'etend donc de mani\`ere usuelle en une 
notion de morphismes entre $et$-champs d\'eriv\'es $(n,li)$-g\'eom\'etriques (voir par exemple
\cite[\S 1.3.6]{hagII}).

Pour d\'emontrer la proposition on consid\`ere un carr\'e commutatif
$$\xymatrix{
Y=\mathbb{R}\underline{Spec}\, B \ar[r]^-{u} \ar[d] & F' \ar[d] \\
X=\mathbb{R}\underline{Spec}\, A \ar[r]_-{v} & F,}$$
avec $B$ une $A$-alg\`ebre simpliciale plate et finie (de rang $m$). Il nous faut montrer que 
le lieu dans $X$ au-dessus du quel le morphisme $u$ est quasi-lisse est un ouvert $U\subset X$. Comme
ceci est une assertion locale pour la topologie \'etale sur $X$, et que localement sur $X_{et}$ 
le morphisme $u$ se factorise par un $(n,li)$-atlas de $F'\times_{F}^{h}X$, on se ram\`ene au cas
o\`u $F=X$ (et $v=id$) et $F'=\mathbb{R}\underline{Spec}\, C$ est affine. On dipose donc d'un diagramme
commutatif d'affines
$$\xymatrix{
Y=\mathbb{R}\underline{Spec}\, B \ar[r]^-{u} \ar[d] & Z=\mathbb{R}\underline{Spec}\, C \ar[dl] \\
X=\mathbb{R}\underline{Spec}\, A, &}$$
avec $B$ plate et finie sur $A$ et $C$ une $A$-alg\`ebre simpliciale plate et presque de pr\'esentation finie. 
On consid\`ere $U\subset X$ le sous-objet de $X$ qui pour $A'\in sk-CAlg$ consiste en 
tous les morphismes $A \longrightarrow A'$ tels que le morphisme induit
$$C\otimes_{A}^{\mathbb{L}}A' \longrightarrow B\otimes_{A}^{\mathbb{L}}A'$$
soit quasi-lisse. On doit montrer que $U$ est un ouvert Zariski de $X$. Pour cela, on utilise le lemme suivant.

\begin{lem}\label{l3}
Soit $A \in sk-CAlg$ et $M\in Ho(sA-Mod)$. Alors $M$ est parfait d'amplitude contenue 
dans $[a,0]$ si et seulement si $M\otimes_{A}^{\mathbb{L}}\pi_{0}(A)$ est parfait et 
d'amplitude contenue dans $[a,0]$ en tant que $\pi_{0}(A)$-module simplicial. 
\end{lem}

\textit{Preuve:} Le n\'ecessit\'e se d\'eduit du fait qu'\^etre parfait d'amplitude donn\'ee est une propri\'et\'e stable par changement de bases. Pour la suffisance on proc\`ede par r\'ecurence sur l'amplitude. 
Pour $a=0$ c'est l'\'enonc\'e \cite[Lem. 2.2.2.2]{hagII}.
Dire que $M\otimes_{A}^{\mathbb{L}}\pi_{0}(A)$ est parfait d'amplitude contenue dans $[a,0]$ \'equivaut \`a dire qu'il 
existe un morphisme 
$$p : \pi_{0}(A)^{n} \longrightarrow M\otimes_{A}^{\mathbb{L}}\pi_{0}(A)$$ 
dont la cofibre homotopique est parfaite et d'amplitude contenue dans $[a,-1]$. Le morphisme $p$ se rel\`eve de fa\c{c}on unique en un morphisme
$$p' : A^{n} \longrightarrow M$$
dont la cofibre est un $A$-module simplicial $N$ tel que $N\otimes_{A}^{\mathbb{L}}\pi_{0}(A)$ est parfait 
et d'amplitude contenue dans $[a,-1]$. Par r\'ecurrence $N$ est donc parfait et d'amplitude contenue dans $[a,-1]$. 
La suite exacte de cofibrations $\xymatrix{A^{n} \ar[r] & M \ar[r] & N}$ implique donc que 
$M$ est parfait et d'amplitude contenue dans $[a,0]$.  \hfill $\Box$ \\

Revenons \`a notre diagramme commutatif
$$\xymatrix{
Y=\mathbb{R}\underline{Spec}\, B \ar[r]^-{u} \ar[d] & Z=\mathbb{R}\underline{Spec}\, C \ar[dl] \\
X=\mathbb{R}\underline{Spec}\, A. &}$$
Le lemme pr\'ec\'edent implique que $u$ est quasi-lisse si et seulement 
si $\mathbb{L}_{u}\otimes_{B}^{\mathbb{L}}\pi_{0}(B)$ est parfait d'amplitude dans $[-1,0]$. 
Comme $B$ et $C$ sont plates sur $A$ on a 
$$\pi_{0}(B)\simeq B\otimes_{A}^{\mathbb{L}}\pi_{0}A \qquad \pi_{0}(C)\simeq C\otimes_{A}^{\mathbb{L}}\pi_{0}A,$$
et donc
$$\mathbb{L}_{u}\otimes_{B}^{\mathbb{L}}\pi_{0}(B)\simeq 
\mathbb{L}_{\pi_{0}(u)},$$
o\`u $\pi_{0}(u) : \pi_{0}(C) \longrightarrow \pi_{0}(B)$ est le morphisme induit. Comme les topologies de Zariski de 
$A$ et de $\pi_{0}(A)$ coincident on voit qu'il suffit de montrer que le lieu dans le sch\'ema $Spec\, \pi_{0}(A)$, au-dessus du quel $\pi_{0}(u)$ est quasi-lisse, est un ouvert Zariksi. En d'autres termes, nous avons ramen\'e le probl\`eme
au cas o\`u $A$, $B$ et $C$ sont des $k$-alg\`ebres commutatives non-simpliciales. En utilisant qu'un morphisme 
de pr\'esentation entre sch\'emas affines est quasi-lisse si et seulement s'il est l.c.i., on voit que l'\'enonc\'e devient alors le fait bien connu suivant, dans le cadre des sch\'emas affines (non-d\'eriv\'es).

\begin{lem}\label{l4}
Soit 
$$\xymatrix{
Y=Spec\, B \ar[r]^-{u} \ar[d]_-{p} & Z=Spec\, C \ar[dl]^-{q} \\
X=Spec\, A, &}$$
un diagramme commutatif de sch\'emas affines,
avec $p$ plat et fini, et $q$ plat et de pr\'esentation finie. Soit $x\in X$ un point tel que 
le morphisme induit sur les fibres $Y_{x} \longrightarrow Z_{x}$ soit de locale intersection compl\`ete. Alors, $u$ est 
de locale intersection compl\`ete au-dessus d'un voisinage Zariski de $p(x) \in X$.
\end{lem}

\textit{Preuve:} Comme tout est de pr\'esentation finie au-dessus de $A$ on se ram\`ene, par un argument standard, au cas 
o\`u tous anneaux en jeu sont noeth\'eriens 
(voir \cite[Cor. 11.2.6.1]{egaIV-3}). Comme $p$ est fini, et donc propre, il suffit de montrer que
si $x\in X$ est tel que $u_{x} : Y_{x} \longrightarrow Z_{x}$ soit l.c.i., il existe un voisinage ouvert $V$ de $Y_{x}$ dans $Y$ tel que 
$u$ soit l.c.i. sur $V$. Le morphisme $u_{x}$ \'etant l.c.i., on peut trouver, localement sur $Y_{x}$, une factorisation
$$u_{x} : \xymatrix{Y_{x} \ar[r]^-{j} & Z_{x}\times \mathbb{A}^{m} \ar[r] & Z_{x}}$$
avec $j$ une immersion ferm\'ee r\'eguli\`ere. L'immersion $j$ se rel\`eve un un morphisme
$i : Y \longrightarrow Z\times \mathbb{A}^{m}$ qui factorise $u$. Soit $(f_{1},\dots,f_{r})$ des g\'en\'erateurs 
de l'id\'eal d\'efinissant $Y$ dans $Z\times \mathbb{A}^{m}$, et formons la dg-alg\`ebre de Koszul 
associ\'ee $K(f)$. On dipose d'une augmentation naturelle $K(f) \longrightarrow B$, que l'on consid\`ere
comme un morphisme de complexes coh\'erents born\'es sur $Z\times \mathbb{A}^{m}$. Par hypoth\`ese ce morphisme de complexe
est un quasi-isomorphisme de complexes lorsqu'il est restreint \`a $Y_{x}$, et reste donc un quasi-isomorphisme
sur un voisinage ouvert de $Y_{x}$ dans $Y$. Mais ceci implique que $u$ est l.c.i. sur un voisinage ouvert 
de $Y_{x}$.  \hfill $\Box$ \\

Cela termine la preuve de la proposition \ref{p6}.
\hfill $\Box$ \\

\subsection{Preuve du th\'eor\`eme}

Nous sommes maintenant en mesure de donner une preuve du th\'eor\`eme \ref{t1}. Pour cela, nous supposons que 
$\phi_{n}$ est une \'equivalence de cat\'egorie, pour un $n\geq 0$. Pour montrer que $\phi_{n+1}$ est 
aussi une \'equivalence il suffit, d'apr\`es \ref{l1}, de montrer que $\phi_{n+1}$ est essentiellement surjectif. \\

Soit donc $F\in dSt_{fppf}^{n+1,pl}(k)$. On sait que le morphisme diagonal 
$F \longrightarrow F\times^{h}F$ est $(n,pl)$-repr\'esentable, et donc aussi $(n,li)$-repr\'esentable par r\'ecurrence.
Soit $\{U_{i}\}$ une famille d'affines et 
$$p=\coprod_{i}p_{i} : U:=\coprod_{i}U_{i} \longrightarrow F$$ 
un $(n+1,pl)$-atlas. Pour tout $i$, et tout entier $m$, on dispose 
du $et$-champ d\'eriv\'e $\underline{QSect}_{p_{i},m}^{str}$ des quasi-sections strictes de de rang $m$ du morphsme $p_{i}$. 
On note $V_{i,m}\subset \underline{QSect}_{p_{i},m}^{str}$
le sous-champ des quasi-sections quasi-lisses. D'apr\`es la proposition \ref{p1}, le morphisme
$V_{i,m} \longrightarrow F$ est $(n+1,li)$-repr\'esentable, car Zariski ouvert dans un $et$-champ d\'eriv\'e
$(n,li)$-repr\'esentable au-dessus de $F$ (le seul cas o\`u il est n\'ecessaire de passer de $n$ \`a $n+1$ est lorsque
$n=0$, l'ouvert $V_{i,m}$ n'\'etant pas forc\'ement affine). On consid\`ere le morphisme de projection
$$V:=\coprod_{i,m} : V_{i,m} \longrightarrow F.$$
Pour terminer la preuve du th\'eor\`eme \ref{t1} il nous suffit de montrer les trois assertions suivantes.
\begin{enumerate}
\item Chaque $et$-champ $V_{i,m}$ est $(n+1,li)$-g\'eom\'etrique.
\item Chaque morphisme $V_{i,m} \longrightarrow F$ est $(n+1,li)$-repr\'esentable et lisse. 
\item Pour tout corps alg\'ebriquement clos $L$, le morphisme induit
$$\coprod_{i,m} : V_{i,m}(L) \longrightarrow F(L)$$
est surjectif sur les composantes connexes. 
\end{enumerate}

Montrer les propri\'et\'es $(1)-(3)$ ci-dessus impliquera que $F$ est $(n+1,li)$-g\'eom\'etrique: un $(n+1,li)$-atlas 
pour $V$ induisant un $(n+1,li)$-atlas pour $F$. \\

$(1)$ Il suffit de montrer que les $et$-champs $\underline{QSect}_{p_{i},m}^{str}$ sont
$(n+1,li)$-g\'eom\'etriques, car les $V_{i,m}$ sont des ouverts des $\underline{QSect}_{p_{i},m}^{str}$.
On sait, d'apr\`es la proposition \ref{p6} que la projection $\underline{QSect}_{p_{i},m}^{str} \longrightarrow 
F$ est $(n,li)$-g\'eom\'etrique. Par construction, il existe un diagramme commutatif de $et$-champs d\'eriv\'es
$$\xymatrix{
\widetilde{\underline{QSect}_{p_{i},m}^{str}} \ar[r]^-{u} \ar[d] & U_{i}\ar[d] \\
\underline{QSect}_{p_{i},m}^{str} \ar[r] & F,}$$
o\`u $\widetilde{\underline{QSect}_{p_{i},m}^{str}} \longrightarrow \underline{QSect}_{p_{i},m}^{str}$
est un morphisme plat et fini (il s'agit de l'image r\'eciproque de l'objet universel
$\widetilde{\underline{Fin}_{m}^{str}} \longrightarrow \underline{Fin}_{m}^{str}$). Comme
$\widetilde{\underline{QSect}_{p_{i},m}^{str}}$ et $U_{i}$ sont $(n,li)$-repr\'esentables au-dessus de $F$, il s'en suit que
le morphisme $u$ est aussi $(n,li)$-repr\'esentable. Ceci implique que 
$\widetilde{\underline{QSect}_{p_{i},m}^{str}}$ est un $et$-champ d\'eriv\'e $(n,li)$-g\'eom\'etrique. Soit 
$Y_{j}$ des affines et 
$$\coprod_{j}Y_{j} \longrightarrow \widetilde{\underline{QSect}_{p_{i},m}^{str}}$$
un $(n,li)$-atlas, en supposant que $n>0$. Le morphisme compos\'e
$$\coprod_{j}Y_{j} \longrightarrow \underline{QSect}_{p_{i},m}^{str}$$
est alors un $(n,pl)$-atlas, ce qui par l'hypoth\`ese de r\'ecurrence implique que 
$\underline{QSect}_{p_{i},m}^{str}$ est $(n,li)$-g\'eom\'etrique si $n>0$. Il faut ici prendre garde au cas 
$n=0$, qui doit se traiter de mani\`ere ind\'ependante. Dans ce cas on consid\`ere le diagramme
homotopiquement cart\'esien suivant
$$\xymatrix{
\widetilde{\underline{QSect}_{p_{i},m}^{str}} \ar[r] \ar[d] & \underline{QSect}_{p_{i},m}^{str} \ar[d]\\
\widetilde{\underline{Fin}_{m}^{str}} \ar[r]_-{f} & \underline{Fin}_{m}^{str},}$$
o\`u $f$ est la famille universelle des morphismes plats finis de rang $m$. On a vu que 
$\widetilde{\underline{QSect}_{p_{i},m}^{str}}$ \'etait affine, ainsi le $et$-champ d\'eriv\'e 
$\underline{QSect}_{p_{i},m}^{str}$ est localement, pour la topologie plate sur $\underline{Fin}_{m}^{str}$, 
affine. On sait, d'apr\`es \cite[Prop. 1.3.2.8]{hagII} que cela implique que $\underline{QSect}_{p_{i},m}^{str}$ est affine. \\

$(2)$ Les morphismes $V_{i,m} \longrightarrow F$ sont $(n+1,li)$-repr\'esentables, d'apr\`es $(1)$ et car la diagonale de 
$F$ est $(n,li)$-repr\'esentable (et donc $(n+1,li)$-repr\'esentable). Par changement de bases sur $F$ on 
se ram\`ene \`a montrer le lemme suivant.

\begin{lem}\label{l6}
Soit $f : F \longrightarrow X$ un morphisme $(n+1,li)$-repr\'esentable plat et presque de pr\'esentation finie avec 
$X$ affine. Le morphisme naturel
$$\underline{QSect}_{f,m}^{str,ql} \longrightarrow X$$
est lisse. 
\end{lem}

\textit{Preuve du lemme:} Nous utiliserons le crit\`ere \cite[Cor. 2.2.5.3]{hagII}. On consid\`ere la factorisation naturelle
$$\xymatrix{\underline{QSect}_{f,m}^{str,ql} \ar[r]^-{p} & X\times^{h}\underline{Fin}_{m}^{str} \ar[r]^-{q} & X,}$$
et les morphismes sur les tronqu\'es associ\'es
$$\xymatrix{t_{0}(\underline{QSect}_{f,m}^{str,ql}) \ar[r]^-{p} & t_{0}(X)\times t_{0}(\underline{Fin}_{m}^{str}) \ar[r]^-{q} & t_{0}(X).}$$
Pour voir que $t_{0}(\underline{QSect}_{f,m}^{str,ql})$ est localement de pr\'esentation finie sur $t_{0}(X)$, il suffit
de voir que $t_{0}(\underline{Fin}_{m}^{str})$ est un sch\'ema affine (non-d\'eriv\'e) de pr\'esentation finie sur $Spec\, k$, et que $p$ est localement de pr\'esentation finie. Le sch\'ema 
$t_{0}(\underline{Fin}_{m}^{str})$  classifie
les structures de $k$-alg\`ebres commutatives sur $k^{m}$, et est donc de pr\'esentation finie. Il reste \`a voir que 
$p$ est un morphisme localement de pr\'esentation finie de $et$-champs non-d\'eriv\'es. 

Pour cela, soit $Y$ un sch\'ema affine et $Y \longrightarrow t_{0}(X)\times t_{0}(\underline{Fin}_{m}^{str})$ un morphisme correspondant \`a un morphisme $Y\rightarrow t_{0}(X)$ et un morphisme $Y' \rightarrow Y$ libre de rang $m$. Le produit fibr\'e homotopique 
$$t_{0}(\underline{QSect}_{f,m}^{str,ql})\times_{t_{0}(X)\times t_{0}(\underline{Fin}_{m}^{str})}^{h}t_{0}(X)$$
est un ouvert du $et$-champ (non-d\'eriv\'e) $\underline{Map}_{/Y}(Y',F\times_{t_{0}(X)}^{h}Y)$, des $Y$-morphismes
de $Y'$ vers $F$ (correspondant au sous-champ des morphismes quasi-lisses). Au cours de la preuve de la proposition \ref{p5} nous avons vu qu'un $(n,li)$-atlas de $\underline{Map}_{/Y}(Y',F\times_{t_{0}(X)}^{h}Y)$ \'etait donn\'e 
par $\underline{Map}_{/Y}(Y',t_{0}(U))$, o\`u $U$ est un $(n,li)$-atlas de $F$. Ainsi, la locale pr\'esentabilit\'e du morphisme
$p$ ci-dessus se d\'eduit du fait que pour tout sch\'ema affine $U \rightarrow Y$, le sch\'ema affine
des morphismes $\underline{Map}_{/Y}(Y',U)$ est localement de pr\'esentation finie sur $Y$. 

Nous venons donc de voir que $t_{0}(\underline{QSect}_{f,m}^{str,ql}) \longrightarrow t_{0}(X)$ est localement de pr\'esentation finie. Il reste \`a montrer que les complexes cotangents du morphisme $\underline{QSect}_{f,m}^{str,ql} \longrightarrow
X$ sont parfaits et d'amplitude contenue dans $[0,\infty[$. 
Pour cela, on utilise la factorisation
$$\xymatrix{\underline{QSect}_{f,m}^{str,ql} \ar[r]^-{p} & X\times^{h}\underline{Fin}_{m} \ar[r]^-{q} & X.}$$
Ainsi, pour tout $B \in sk-CAlg$, $Y:=\mathbb{R}\underline{Spec}\, B$, et tout morphisme $v : Y \rightarrow \underline{QSect}_{f,m}^{str,ql}$, on dispose d'un triangle distingu\'e de $B$-modules stables
$$\xymatrix{\mathbb{L}_{X\times^{h}\underline{Fin}_{m}/X,pv} \ar[r] & 
\mathbb{L}_{\underline{QSect}_{f,m}^{str,ql}/X,v} \ar[r] & 
\mathbb{L}_{\underline{QSect}_{f,m}^{str,ql}/X\times^{h}\underline{Fin}_{m},v}.}$$
Supposons que le morphisme $v$ corresponde \`a un diagramme commutatif
$$\xymatrix{
Y' \ar[r]^-{u} \ar[d] & F \ar[d] \\
Y \ar[r] & X,}$$
avec $Y'=\mathbb{R}\underline{Spec}\, B'$ plat et fini sur $Y$, et $u$ quasi-lisse. 

Alors, on a 
$$\mathbb{L}_{\underline{QSect}_{f,m}^{str,ql}/X\times^{h}\underline{Fin}_{m},v}\simeq
\mathbb{L}_{\underline{Map}_{/Y}(Y',F\times_{X}^{h}Y),u}.$$ 
Tout comme pour la proposition \ref{p3} on voit qu'il existe un isomorphisme naturel
$$\mathbb{L}_{\underline{Map}_{/Y}(Y',F\times_{X}^{h}Y),u}
\simeq \mathbb{L}_{F,u}\otimes_{B}^{\mathbb{L}}(B')^{\vee} \in Ho(Sp(sB-Mod)).$$
On a de m\^eme, 
$$\mathbb{L}_{X\times^{h}\underline{Fin}_{m}^{str}/X,pv}\simeq \mathbb{L}_{\underline{Fin}_{m},w} \in Ho(Sp(sB-Mod)),$$
o\`u $w : Y \longrightarrow \underline{Fin}_{m}$ est le morphisme induit, correspondant au morphisme fini et plat 
$Y' \rightarrow Y$. Ainsi, d'apr\`es la proposition \ref{p3} 
$\mathbb{L}_{\underline{QSect}_{f,m}^{str,ql}/X,v}$ est la fibre homotopique, dans $Ho(Sp(sB-Mod))$, du morphisme naturel
$$\mathbb{L}_{F,u}\otimes_{B'}^{\mathbb{L}}(B')^{\vee}\longrightarrow  \mathbb{L}_{B'/B}\otimes_{B'}^{\mathbb{L}}(B')^{\vee},$$
(induit par le morphisme $\mathbb{L}_{F,u} \longrightarrow \mathbb{L}_{B'/B}$ induit par $u$). 
Comme $u$ est quasi-lisse cette fibre homotopique est parfaite et d'amplitude contenue dans $[0,\infty[$, en tant que 
$B'$-module. Mais comme $B'$ est projective et de rang fini sur $B$ le $B$-module
$\mathbb{L}_{\underline{QSect}_{f,m}^{str,ql}/X,v}$ est aussi parfait et d'amplitude contenue dans 
$[0,\infty[$. \hfill $\Box$ \\

$(3)$ Soit $x : Spec\, L \longrightarrow F$ un point g\'eom\'etrique de $F$, et notons 
$$\coprod_{i}q_{i} : \coprod_{i}G_{i} \longrightarrow Spec\, L$$ 
le changement de base du $(n+1,pl)$-atlas $\{U_{i}\}$ le long de $x$. Remarquons que chaque $G_{i}$
est plat sur $Spec\, L$, et donc est un $et$-champ non-d\'eriv\'e (i.e. \'equivalent \`a son tronqu\'e).
Il nous suffit de montrer que 
$\underline{QSec}_{q_{i},m}^{str,ql}(L)$ est non-vide pour un $i$ et un $m$. Pour cela, soit 
$V_{i}$ un affine non-vide et $V_{i} \longrightarrow G_{i}$ un morphisme lisse (pour un $i$ fix\'e quelconque).
Comme $V_{i}$ est plat sur $Spec\, L$, il s'agit d'un sch\'ema affine non-d\'eriv\'e
Soit alors $y \in V_{i}$ un point Cohen-MacCauley de $V_{i}$, et $(f_{1},\dots,f_{r})$
une suite r\'eguli\`ere maximale en $y$. La $L$-alg\`ebre $R:=\mathcal{O}_{V_{i},y}/(f_{1},\dots,f_{r})$
est finie sur $L$, de dimension un entier $m$, et part d\'efinition le morphisme naturel 
$Spec\, R \longrightarrow V_{i}$ est l.c.i. Ces donn\'ees fournissent une quasi-section quasi-lisse 
du morphisme $V_{i} \longrightarrow Spec\, L$. Comme le morphisme $V_{i} \longrightarrow G_{i}$ est lisse, 
l'image de cette quasi-section dans $G_{i}$ fournit un \'el\'ement dans 
$\pi_{0}(\underline{QSec}_{q_{i},m}^{str,ql}(L))$, montrant ainsi que $\underline{QSec}_{q_{i},m}^{str,ql}(L)$
n'est pas vide. \\

Nous venons de voir que les assertions $(1)-(3)$ \'etaient satisfaites, ce qui finit de d\'emontrer le th\'eor\`eme
\ref{t1}. \\

\section{Application \`a la comparaison entre cohomologies \'etales et plates}

Pour un $fppf$-champ d\'eriv\'e, nous savons maintenant qu'\^etre $(n,pl)$-g\'eom\'etrique et 
$(n,li)$-g\'eom\'etrique sont deux conditions \'equivalentes. Nous dirons alors simplement 
\emph{\^etre $n$-g\'eom\'etrique}. Nous dirons aussi \emph{\^etre g\'eom\'etrique} pour signifier 
\emph{\^etre $n$-g\'eom\'etrique pour un entier $n$}. \\

Rappelons que pour une cat\'egorie de mod\`eles $M$ (ou plus g\'en\'eralement pour une sous-cat\'egorie pleine d'une cat\'egorie de mod\`eles, stable par \'equivalence) on dispose d'une notion 
d'objet en groupo\"\i des de Segal $X_{*}$ dans $M$ (voir \cite[Def. 4.9.1]{hagI}). 
Nous dirons qu'un 
objet en groupo\"\i des de Segal $X_{*}$ est \emph{un objet en groupes dans $M$} si l'objet 
$X_{0}\in M$ est \'equivalent \`a l'objet final $*$.
Nous appelerons alors \textit{$\tau$-champ d\'eriv\'e en groupes} un objet en groupes
$X_{*}$ dans la cat\'egorie de mod\`eles $dAff_{k}^{\sim,\tau}$, tel que chaque $X_{n}$ soit un $\tau$-champ. 
Nous abuserons souvent du fait de d\'esigner l'objet 
$X_{*}$ par son objet sous-jacent $G:=X_{1} \in dAff_{k}^{\sim,\tau}$. 
Pour un $\tau$-champ d\'eriv\'e en groupes $G$ on dispose de son $\tau$-champ d\'eriv\'e classifiant 
$$K_{\tau}(G,1):=Hocolim_{[n]\in \Delta^{op}}X_{n} \in dSt_{\tau}(k),$$ 
(voir \cite[Def. 4.9.1]{hagI}).

Par d\'efinition, un \emph{$\tau$-champ d\'eriv\'e en $m$-groupes} est un objet en groupes dans 
la cat\'egorie de (pseudo, voir \cite[Def. 4.1.1]{hagI}) mod\`eles des $\tau$-champs d\'eriv\'es en $(m-1)$-groupes. Pour un tel objet 
$G$, son $\tau$-champ d\'eriv\'e classifiant $K_{\tau}(G,1)$ poss\`ede une structure induite 
de $\tau$-champ d\'eriv\'e en $(m-1)$-groupes. Par it\'eration on obtient ainsi un $\tau$-champ classifiant 
$K_{\tau}(G,m)$.

\begin{df}\label{d3}
Soit $X$ un $\tau$-champ d\'eriv\'e et $G$ un $\tau$-champ d\'eriv\'e en $m$-groupes. La cohomologie de $X$ \`a coefficients dans $G$ est d\'efinie par 
$$H^{m-i}_{\tau}(X,G):=\pi_{i}(Map(X,K_{\tau}(G,m))).$$
\end{df}

Tout d'abord, $K_{\tau}(G,m)$ ne poss\`edant aucune structure de groupe induite les $H^{j}_{\tau}(X,G)$ sont 
des groupes ab\'eliens pour $j\leq m-2$, $H^{m-1}_{\tau}(X,G)$ est un groupe, et 
$H^{m}_{\tau}(X,G)$ est un ensemble point\'e. 
On remarque que $H^{j}_{\tau}(X,G)$ n'est d\'efinie que pour $j\leq m$, mais peut \^etre non-nul pour 
$j<0$ (mais ce qui n'arrive que lorsque $X$ n'est pas un champ tronqu\'e). 
Enfin, lorsque $G$ est un $k$-sch\'ema en groupes (resp. en groupes ab\'eliens) et $X$ un $k$-sch\'ema, 
les $H_{\tau}^{j}(X,G)$ coincident avec la cohomologie de $X$ \`a coefficients dans 
$G$ pour la topologie $\tau$ au sens usuel de la cohomologie des sch\'emas. \\

Une cons\'equence du th\'eor\`eme \ref{t1} est le corollaire suivant. Il affirme en particulier que la cohomologie
$fppf$ et \'etale, \`a coefficients dans un sch\'ema en groupes plats de pr\'esentation fini, ne diff\`erent 
essentiellement qu'en degr\'e $1$. 

\begin{cor}\label{c1}
Soit $G$ un $fppf$-champ d\'eriv\'e en $m$-groupes.
\begin{enumerate}
\item Si $G$ est $n$-g\'eom\'etrique, plat et localement de pr\'esentation presque finie sur $Spec\, k$, alors
$K_{fppf}(G,m)$ est $(n+m)$-g\'eom\'etrique plat et de pr\'esentation preque finie sur $Spec\, k$. Si $m>0$, alors
$K_{fppf}(G,m)$ est lisse sur $Spec\, k$. 
\item Si $G$ est g\'eom\'etrique et plat de pr\'esentation presque finie sur $Spec\, k$, alors
le morphisme $G \longrightarrow Spec\, k$ est quasi-lisse.
\item Si $G$ est g\'eom\'etrique et lisse sur $Spec\, k$, alors pour tout $k$-sch\'ema $X$, le morphisme naturel
$$H^{i}_{et}(X,G) \longrightarrow H_{fppf}^{i}(X,G)$$
est bijectif pour tout $i\leq m$. 
\item Supposons que $G$ soit un $k$-sch\'emas en groupes ab\'eliens (ou plus g\'en\'eralement un 
$k$-espace alg\'ebrique en groupes ab\'eliens), plat et localement de pr\'esentation sur $Spec\, k$.
Alors, pour tout $k$-sch\'ema $X$ il existe une suite exacte longue fonctorielle en $X$ et en $G$
$$\xymatrix{
H^{i-2}_{et}(X,\underline{H}^{1}_{fppf}(-,G)) \ar[r] & H^{i}_{et}(X,G) \ar[r] & H^{i}_{fppf}(X,G) \ar[r] & H^{i-1}_{et}(X,\underline{H}^{1}_{fppf}(-,G)) \ar[r] & 
H^{i+1}_{et}(X,G), }$$
o\`u $\underline{H}^{1}_{fppf}(-,G)$ est le faisceau, pour la topologie \'etale, associ\'e 
au pr\'efaisceau $U \mapsto H^{1}_{fppf}(U,G)$.
\item Supposons que $G$ soit un $k$-sch\'emas en groupes ab\'eliens (ou plus g\'en\'eralement un 
$k$-espace alg\'ebrique en groupes ab\'eliens), plat et localement de pr\'esentation sur $Spec\, k$.
Alors, pour tout $k$-sch\'ema $X$ et toute classe $\alpha \in H^{i}_{fppf}(X,G)$, avec $i>1$, il existe 
un recouvrement \'etale $u : X' \rightarrow X$ tel que $u^{*}(\alpha)=0$. En particulier, si $X$ est 
hens\'elien strict alors $H^{i}_{fppf}(X,G)=0$ pour tout $i>1$. 
\end{enumerate}
\end{cor}

\textit{Preuve:} $(1)$ Le $fppf$-champ d\'eriv\'e $K_{fppf}(G,m)$ est, par des applications successives
de \cite[Prop. 1.3.4.2]{hagII}, $(m+n,pl)$-g\'eom\'etrique, et donc $(n+m)$-g\'eom\'etrique. Pour tout $m>0$ le point global naturel
$Spec\, k \longrightarrow K_{fppf}(G,m)$ est plat et localement de pr\'esentation presque finie. 
Ainsi, $K_{fppf}(G,m)$ est localement pour la topologie $fppf$ lisse sur $Spec\, k$, ce qui d'apr\`es 
la proposition \ref{p1} implique qu'il est lisse sur $Spec\, k$. \\

$(2)$ D'apr\`es le point pr\'ec\'edent $K_{fppf}(G,1)$ est lisse. Le $et$-champ d\'eriv\'e
$G$, qui s'\'ecrit aussi $*\times_{K_{fppf}(G,1)}^{h}*$, est donc un produit fibr\'e de 
$et$-champs d\'eriv\'es lisses. Il est donc quasi-lisse. \\

$(3)$ Le $et$-champ d\'eriv\'e $K_{et}(G,m)$ est, par une utilisation r\'ep\'et\'ee 
de \cite[Prop. 1.3.4.2]{hagII}, $(n+m,li)$-g\'eom\'etrique. Le lemme \ref{l1} $(1)$ implique que pour $et$-champ 
$X$ (et donc, en particulier, pour 
tout $k$-sch\'ema $X$), le morphisme naturel 
$$H^{i}_{et}(X,G)=[X,K_{et}(G,i)] \longrightarrow [X,K_{fppf}(G,i)]=H_{fppf}^{i}(X,G)$$
est bijectif, pour tout $i\leq m$. \\

$(4)$ On consid\`ere $K_{fppf}(G,1)$. D'apr\`es le point $(1)$ c'est un champ d\'eriv\'e g\'eom\'etrique en 
$m$-groupes (pour tout $m$ car $G$ est ab\'elien), lisse sur $Spec\, k$. Soit $m$ un entier assez grand.
Le $et$-champ d\'eriv\'e
$K_{et}(K_{fppf}(G,1),m-1)$ est donc g\'eom\'etrique, et se trouve donc 
\^etre un $fppf$-champ (voir le lemme \ref{l1}). On a donc
$$K_{fppf}(G,m)\simeq K_{fppf}(K_{fppf}(G,1),m-1)\simeq
K_{et}(K_{fppf}(G,1),m-1).$$
Ainsi, les faisceaux d'homotopie, pour la topologie \'etale,  du $et$-champ tronqu\'e $t_{0}(K_{fppf}(G,m))$ sont donn\'es par
$$\pi_{m}(t_{0}(K_{fppf}(G,m)))\simeq G \qquad \pi_{m-1}(t_{0}(K_{fppf}(G,m)))\simeq \underline{H}^{1}_{fppf}(-,G) \qquad
\pi_{m-1}(t_{0}(K_{fppf}(G,m)))=0 \; pour \; i<m-1.$$
Il existe donc une suite exacte de fibration de $et$-champs non-d\'eriv\'es
$$t_{0}(K_{et}(G,m)) \longrightarrow t_{0}(K_{fppf}(G,m)) \longrightarrow t_{0}(K_{et}(\underline{H}^{1}_{fppf}(-,G),m-1)).$$
Pour tout $k$-sch\'ema, cette suite excate de fibration induit une suite exacte de fibration
d'ensembles simpliciaux
$$Map(X,K_{et}(G,m)) \longrightarrow Map(X,K_{fppf}(G,m)) \longrightarrow Map(X,K_{et}(\underline{H}^{1}_{fppf}(-,G),m-1)),$$
obtenue en remarquant que comme $X$ est tronqu\'e, on a $Map(X,F)\simeq Map(X,t_{0}(F))$
pour tout $et$-champ d\'eriv\'e $F$. La suite exacte longue en homotopie, associ\'ee \`a cette suite 
exacte de fibration, est la suite exacte cherch\'ee.  \\

$(5)$ Se d\'eduit, par exemple, du point $(4)$ (ou encore de la derni\`ere partie du point $(1)$). \hfill $\Box$ \\

\begin{rmk}
\begin{enumerate}
\item
\emph{Comme nous l'avons vu au cours de la preuve, le point $(3)$ du corollaire pr\'ec\'edent est vrai pour tout $et$-champ
d\'eriv\'e $X$. Cela est \'egalement le cas du point $(5)$. De m\^eme, le point $(4)$ est vrai pour tout 
$et$-champ tronqu\'e $X$}.
\item
\emph{Le point $(4)$ peut aussi s'exprimer sous la forme $\mathbb{R}^{i}f_{*}(G)=0$, pour
tout $i>1$, et tout sch\'ema en groupes ab\'eliens plats localement de pr\'esentation finie $G$, et o\`u 
$f : Aff_{k}^{\sim,fppf} \longrightarrow Aff_{k}^{\sim,et}$ est le morphisme
g\'eom\'etrique de passage de la topologie $fppf$ \`a la topologie \'etale.}
\end{enumerate}
\end{rmk}

Pour terminer, signalons aussi le corollaire suivant. 

\begin{cor}\label{c2}
Soit $A$ un anneau local hens\'elien de corps r\'esiduel $k$. Supposons que $A$ soit un anneau excellent. 
Soit $G$ un espace alg\'ebrique en groupes ab\'eliens, plat et localement de pr\'esentation finie sur $Spec\, A$, de fibre sp\'eciale $G_{0}:=G\otimes_{A}k$. Alors, le morphisme de restriction
$$H^{i}_{fppf}(Spec\, A,G) \longrightarrow H^{i}_{fppf}(Spec\, k,G_{0})$$
est surjectif pour $i=1$, et un isomorphisme pour tout $i>1$. Si le groupe $G$ est de plus lisse sur $Spec\, A$, alors ce morphisme est aussi un isomorphisme pour $i=1$ et surjectif pour $i=0$.  
\end{cor}

\textit{Preuve:} On travaille dans $St(A)$, la cat\'egorie des champs (disons $fppf$) sur $Spec\, A$. 
Commen\c{c}ons par montrer que les morphismes en question sont surjectifs. Pour cela, on consid\`ere le $i$-champ $K_{fppf}(G,i)$, pour un $i\geq 1$. On sait, d'apr\`es le corollaire pr\'ec\'edent que c'est un champ g\'eom\'etrique et lisse sur $Spec\, A$. Ainsi, pour montrer que le morphisme induit
$$H^{i}_{fppf}(Spec\, A,G)\simeq [Spec\, A,K_{fppf}(G,i)] \longrightarrow [Spec\, k,K_{fppf}(G,i)]\simeq H^{i}_{fppf}(Spec\, k,G)$$
est surjectif il suffit de montrer le lemme plus g\'en\'eral suivant. 

\begin{lem}\label{lc2}
Soit $F$ un ($i$-)champ g\'eom\'etrique et lisse sur $Spec\, A$. Alors le morphisme naturel
$$F(A) \longrightarrow F(k)$$
induit une application surjective sur les ensembles de composantes connexes. 
\end{lem}

\textit{Preuve du lemme:} 
Soit $\widehat{A}$ le compl\'et\'e de $A$ le long de $A \rightarrow k$, et fixons $x\in F(k)$ un $k$-point.
Comme le champ $F$ est g\'eom\'etrique le morphisme naturel 
$$F(\widehat{A}) \longrightarrow Holim_{k}(F(A_{k}))$$
est une \'equivalence (on note $A_{k}=A/m^{k}$, o\`u $m$ est l'id\'eal maximal de $A$). 
On commence par remarquer, par r\'ecurrence sur $k$, que le morphisme
$$F(A_{k}) \longrightarrow F(A_{k-1})$$
est surjectif sur les composantes connexes. En effet, l'obstruction au fait que la fibre homotopique 
d'un \'el\'ement $x\in F(A_{k-1})$ soit non-vide vit dans le groupe
$$Ext^{1}_{A_{k-1}}(\mathbb{L},m^{k-1}/m^{k}),$$
o\`u $\mathbb{L}$ est le complexe cotangent du morphisme $F \longrightarrow Spec\, A$ pris au point 
$x : Spec\, A_{k-1} \longrightarrow F$. Or, comme ce morphisme est lisse, $\mathbb{L}$ est d'amplitude 
positive, et donc $Ext^{1}_{A_{k-1}}(\mathbb{L},m^{k-1}/m^{k})=0$ (voir par exemple 
\cite[1.4.2, 2.2.5]{hagII}). Ceci montre que toutes les fibres homotopiques de 
$$F(A_{k}) \longrightarrow F(A_{k-1})$$
sont non vides, et donc que ce morphisme est surjectif sur les composantes connexes. En passant \`a la limite sur $k$, on trouve donc que
$$F(\widehat{A}) \longrightarrow F(k)$$
est aussi surjectif sur les composantes connexes. 

Soit maintenant $x\in F(k)$, et choisissons $\widehat{x} \in F(\widehat{A})$ un relev\'e.
Comme $A$ est excellent, le morphisme $A \longrightarrow \widehat{A}$ est r\'egulier, et il s'\'ecrit donc 
comme une colimite filtrante de morphismes lisses $A \longrightarrow B_{\alpha}$. Comme $A$ est hens\'elien, 
il existe, pour tout $\alpha$, une r\'etraction $r_{\alpha} : B_{\alpha} \longrightarrow A$, au-dessus de $k$. De plus, comme $F$ est localement de pr\'esentation finie sur $Spec\, A$, il existe un $\alpha$ assez grand
tel que $\widehat{x}$ se factorise par un point $x_{\alpha}\in F(B_{\alpha})$ au-dessus de $x$. En composant avec la r\'etraction $r_{\alpha}$ on trouve $x'\in F(A)$, qui est un rel\`evement de $x$. Ceci termine la preuve du lemme.
\hfill $\Box$ \\

Le lemme implique donc que 
$$H^{i}_{fppf}(Spec\, A,G) \longrightarrow H^{i}_{fppf}(Spec\, k,G_{0})$$
est surjectif pour $i>0$. Supposons maintenant que $i>1$. Soit $x,y \in H^{i}_{fppf}(Spec\, A,G)$ avec une m\^eme image dans $H^{i}_{fppf}(Spec\, k,G_{0})$. On repr\'esente $x$ et $y$ par deux morphismes
$$x,y : Spec\, A \longrightarrow K_{fppf}(G,i),$$
et on consid\`ere alors le champ sur $Spec\, A$ des \'equivalences entre $x$ et $y$
$$Eq(x,y):=Spec\, A \times_{K_{fppf}(G,i)}^{h}Spec\, A.$$
Le champ $Eq(x,y)$ est un champ g\'eom\'etrique sur $Spec\, A$. Il est de plus localement \'equivalent, pour la topologie $fppf$, \`a $K_{fppf}(G,i-1)$. Comme $i>1$, le champ 
$K_{fppf}(G,i-1)$ est lisse, ce qui montre que 
$Eq(x,y)$ est un champ g\'eom\'etrique et lisse sur $Spec\, A$. De plus, par hypoth\`ese $Eq(x,y)(k)$ est non vide, et le lemme \ref{lc2} implique alors que $Eq(x,y)(A)$ est aussi non vide. En d'autres termes les deux morphismes
$$x,y : Spec\, A \longrightarrow K_{fppf}(G,i)$$
sont \'egaux dans la cat\'egorie homotopique des champs, et donc $x=y$ dans $H^{i}_{fppf}(Spec\, A,G)$. 

Pour terminer, le m\^eme argument que pr\'ec\'edemment montre que $H^{1}_{fppf}(Spec\, A,G) \longrightarrow H^{1}_{fppf}(Spec\, k,G_{0})$ est injectif lorsque $G$ est lisse (car on peut l'appliquer au cas $i=1$). \hfill $\Box$ \\

\end{document}